\documentclass[final,1p,times]{elsarticle}
\journal{E}
\usepackage{tikz-cd}
\usepackage{extpfeil}
\usepackage{amssymb}
\usepackage{amsthm}
\usepackage{latexsym}
\usepackage{amsmath}
\usepackage{color}
\usepackage{graphicx}
\usepackage{indentfirst}
\usepackage{mathrsfs}
\usepackage{lipsum}
\usepackage{float}
\allowdisplaybreaks

\newtheorem{theorem}{\color{black}\indent \textbf{Theorem}}[section]

\newtheorem{proposition}{\color{black}\indent Proposition}[section]
\newtheorem{definition}{\color{black}\indent Definition}[section]
\newtheorem{remark}{\color{black}\indent Remark}[section]
\newtheorem{corollary}{\color{black}\indent Corollary}[section]
\newtheorem{example}{\color{black}\indent Example}[section]

\begin{document}
\begin{frontmatter}

\title{A geometric approach to the (generalized) Noether theorem for Hamiltonian systems on (non)-uniform $q$-contact manifolds}

\author[a,b]{Manuel de León\corref{cor1}}
\author[c]{Cristina Sardon}
\author[d]{Xuefeng Zhao}
\address[a]{Instituto de Ciencias Matematicas, Campus Cantoblanco, Consejo Superior de Investigaciones 
Científicas, C/ Nicolás Cabrera, 13–15, 28049, Madrid, Spain}
\address[b]{Real Academia Española de las Ciencias. C/ Valverde, 22, 28004, Madrid, Spain}
\address[c]{Departamento de Matemática Aplicada, Universidad Politécnica de Madrid, Escuela de Edificación.  Av. Juan de Herrera 6, 28040 Madrid, Spain}
\address[d]{College of Mathematics, Jilin University, Changchun 130012, P. R. China}
\cortext[cor1]{Corresponding author}

\begin{abstract}
In this paper, we develop a unified Hamiltonian framework for dynamical systems with multiple independent dissipation channels based on $q$-contact geometry. For uniform $q$-contact manifolds (satisfying $d\lambda_1 = \cdots = d\lambda_q$), we prove a Liouville-type theorem characterizing the evolution of the natural volume form along the $q$-contact Hamiltonian flow, establish the $q$-contact Noether theorem, and discuss $q$-contact transformations together with their generating functions. For time-dependent systems, we introduce generalized Noether symmetries on the extended phase space and prove a generalized Noether theorem. We then extend the theory to non-uniform $q$-contact structures, where the exterior derivatives of the contact forms are allowed to differ, yielding a family of distinct symplectic structures on the horizontal distribution. By fixing a reference symplectic form $\omega = d\lambda_1|_\xi$, we introduce structure endomorphisms $B_i$ encoding anisotropic dissipation effects, and define channel Hamiltonian vector fields $X_H^i$ for each geometric channel. Under the admissibility condition that the averaged $2$-form $\Omega = \frac{1}{q}\sum_{i=1}^q d\lambda_i$ be non-degenerate on the horizontal distribution, we construct a unique effective Hamiltonian vector field $X_H$, establish the associated non-uniform $q$-contact bracket and the corresponding Noether theorem, derive the evolution formula for the volume form along the effective flow, and introduce the effective structure endomorphism $B = \frac{1}{q}\sum_{i=1}^q B_i$ which captures the averaged geometric influence of the multiple channels. As an application, we embed an elastoplastic damage model with two internal variables into the non-uniform $2$-contact framework, explicitly construct the underlying geometric structure, derive the effective equations of motion, and validate the theoretical predictions through numerical simulations, which confirm the exponential dissipation law and the associated rescaled conserved quantity. 
\end{abstract}

\end{frontmatter}
\vspace{1em}

{\large\bf\raggedright
    Keywords:} $q$-contact manifold, Hamiltonian system, symmetries, Noether theorem, $q$-contact bracket
\medskip

{\large\bf\raggedright
MSC2020 codes}: 53C15, 37J06, 53D99

\section{Introduction}
The geometric description of mechanical systems has a long and distinguished history, beginning with the symplectic formulation of Hamiltonian mechanics, which provides a natural setting for conservative dynamics. On a symplectic manifold $(M,\omega)$, the dynamics of a system with Hamiltonian function $H$ is governed by the Hamilton equations $i_{X_H}\omega = dH$, and the flow preserves both the symplectic form and the Hamiltonian function itself. This framework has proven extraordinarily successful for closed systems, where energy is conserved and the phase space volume is preserved by Liouville's theorem.

However, many systems of physical interest are inherently dissipative. Phenomena such as friction, viscosity, plastic deformation, damping, and external forcing are ubiquitous in real-world applications, ranging from classical mechanics and engineering to geophysics and biophysics. These effects cannot be adequately described within the symplectic paradigm, as they involve energy exchange with the environment and irreversible entropy production. The development of geometric structures capable of incorporating dissipation has therefore become an important theme in mathematical physics.

Contact geometry \cite{Bravetti,Kholodenko,Liu} offers a natural extension of symplectic geometry to dissipative settings. A contact manifold $(M,\eta)$ is an odd-dimensional manifold endowed with a maximally non-integrable hyperplane distribution $\ker\eta$, equivalently characterized by the condition $\eta\wedge(d\eta)^n \neq 0$. The Reeb vector field $R$, uniquely defined by $i_R\eta = 1$ and $i_Rd\eta = 0$, generates the direction of energy dissipation. On a contact manifold, the Hamiltonian vector field $X_H$ associated with a function $H$ is defined by $i_{X_H}\eta = -H$ and $i_{X_H}d\eta = dH - dH(R)\eta$. Unlike the symplectic case, the Hamiltonian function is generally not conserved along the flow; instead, it satisfies $\dot{H} = -H\,dH(R)$, reflecting the dissipative nature of the dynamics. Contact Hamiltonian systems have been extensively studied and applied to various dissipative phenomena, including friction, thermostatistics, quantum systems \cite{Ciaglia2018,Herczeg,Hooft} and general relativity \cite{Cariglia}.

A further generalization, known as $q$-contact geometry \cite{Finamore,Finamore2,Leok,Zhao}, has been introduced to accommodate systems with multiple independent dissipation channels. A $q$-contact manifold $(M,\vec{\lambda},\mathcal{R}\oplus\xi)$ is a $(2n+q)$-dimensional manifold equipped with a collection of $q$ linearly independent $1$-forms $\lambda_1,\ldots,\lambda_q$ such that the horizontal distribution $\xi = \cap_i\ker\lambda_i$ has rank $2n$, each $d\lambda_i|_\xi$ is non-degenerate, and $\ker d\lambda_i = \mathcal{R}$, where $\mathcal{R}$ is the Reeb distribution spanned by the Reeb vector fields $R_1,\ldots,R_q$ satisfying $\lambda_i(R_j)=\delta_{ij}$. This structure naturally models systems with several dissipative mechanisms, multiple constraints, or coupled energy exchanges. The Reeb vector fields generate distinct dissipation channels, and the interplay among them gives rise to a rich geometric structure that generalizes both symplectic and single-contact geometries.

A central theme in mechanics is the relationship between symmetries and conservation laws, famously encapsulated in Noether's theorem \cite{Aldaya,Aldaya2,de4,Zhao1}. In its classical formulation, the theorem states that to every continuous symmetry of a Lagrangian or Hamiltonian system there corresponds a conserved quantity \cite{Cantrijn,Carinena,Carinena2,de2,Prince1,Prince2,Sarlet}. This result has been extended and refined in numerous directions, including time-dependent systems \cite{de0,Prieto,Prieto2}, gauge symmetries \cite{Prinz}, and higher-order Lagrangians \cite{Gracia,Gracia2}. In the context of dissipative systems, however, the notion of conservation must be modified. For contact Hamiltonian systems, it has been shown that symmetries give rise not to conserved quantities but to dissipated quantities, which decay at a rate determined by the Reeb vector field. This observation motivates the development of a Noether-type theorem for $q$-contact geometry, where multiple dissipation channels may contribute to the decay of the associated quantities.

In this paper, we develop a comprehensive Hamiltonian framework for both uniform and non-uniform $q$-contact manifolds. The uniform case, characterized by $d\lambda_1 = \cdots = d\lambda_q$, provides a natural starting point in which all dissipation channels share the same symplectic structure on the horizontal distribution. In this setting, we derive the $q$-contact Hamiltonian equations in local coordinates, establish a Liouville-type theorem that characterizes the evolution of the natural volume form $\mu = \lambda_1\wedge\cdots\wedge\lambda_q\wedge(d\lambda_1)^n$ along the Hamiltonian flow, and prove a Noether theorem that establishes an equivalence between Noether symmetries and dissipated quantities. We also discuss $q$-contact transformations and their generating functions, providing the analogue of canonical transformations in this generalized setting. Furthermore, for time-dependent systems, we introduce generalized Noether symmetries on the extended phase space and prove a corresponding generalized Noether theorem.

We then extend the theory to non-uniform $q$-contact structures, where the exterior derivatives $d\lambda_1,\ldots,d\lambda_q$ are allowed to differ. This relaxation of the uniformity condition is essential for applications in which different dissipation channels exhibit distinct geometric characteristics, such as anisotropic materials or systems with coupled internal variables. In the non-uniform setting, the horizontal distribution carries a family of distinct symplectic structures. By fixing a reference symplectic form $\omega = d\lambda_1|_\xi$, we introduce structure endomorphisms $B_i:\xi\to\xi$, defined by $d\lambda_i(X,Y)=\omega(B_iX,Y)$, which encode the anisotropic dissipation effects of each channel. For each geometric channel $i$, we define a corresponding Hamiltonian vector field $X_H^i$ determined by $i_{X_H^i}d\lambda_i = dH - \sum_j dH(R_j)\lambda_j$ and $\lambda_j(X_H^i)=-H$ for $j=1,\ldots,q$. These channel vector fields share the same Reeb component but differ in their horizontal components, reflecting the distinct symplectic structures associated with each contact form. Under the admissibility condition that the averaged $2$-form $\Omega = \frac{1}{q}\sum_{i=1}^q d\lambda_i$ be non-degenerate on the horizontal distribution, we construct a unique effective Hamiltonian vector field $X_H$ that unifies all geometric channels into a single dynamical system. We establish the associated non-uniform $q$-contact bracket, prove the corresponding Noether theorem, derive the evolution formula for the volume form along the effective flow, and introduce the effective structure endomorphism $B = \frac{1}{q}\sum_{i=1}^q B_i$, which captures the averaged geometric influence of the multiple channels. In the uniform limit, where $B_i = \mathrm{Id}_\xi$ for all $i$, the effective theory reduces to the uniform $q$-contact Hamiltonian framework.

To demonstrate the applicability of the theoretical framework, we construct a detailed physical model: an elastoplastic material with damage, described by two internal variables representing distinct dissipative mechanisms. The phase space is $M = \mathbb{R}^6(\epsilon, p, \alpha_1, \alpha_2, s_1, s_2)$, where $\epsilon$ is the macroscopic strain, $p$ its conjugate momentum, $\alpha_1,\alpha_2$ are internal variables associated with inelastic deformation, and $s_1,s_2$ are dissipative variables corresponding to the two dissipation channels. The Hamiltonian function is given by
\[
H = \frac{p^2}{2m} + \frac{k}{2}(\epsilon-\alpha_1-\alpha_2)^2 + \frac{k_1}{2}\alpha_1^2 + \frac{k_2}{2}\alpha_2^2 + \kappa\alpha_1\alpha_2 + \mu_1 s_1 + \mu_2 s_2,
\]
which includes elastic energy, internal variable coupling, and linear dissipation terms. We explicitly construct a non-uniform $2$-contact structure on this phase space, verify the admissibility condition, derive the effective equations of motion, and perform numerical simulations. The numerical results confirm the theoretical predictions: the Hamiltonian decays exponentially according to $\dot{H} = -(\mu_1+\mu_2)H$, while the rescaled quantity $I(t)=H(t)e^{(\mu_1+\mu_2)t}$ remains conserved up to numerical precision. This example illustrates the coexistence of dissipated and conserved quantities within the non-uniform $q$-contact framework and validates the geometric theory.

The main contributions of this work are threefold. First, we provide a complete Hamiltonian formulation of dissipative systems on uniform $q$-contact manifolds, including Liouville-type and Noether-type theorems, as well as a treatment of time-dependent symmetries. Second, we develop the theory of non-uniform $q$-contact structures, introducing structure endomorphisms, channel Hamiltonian vector fields, and the effective dynamics under the admissibility condition, thereby extending the geometric framework to systems with anisotropic or coupled dissipation channels. Third, we demonstrate the practical utility of the theory through a concrete application to elastoplastic damage modeling, supported by numerical simulations that validate the predicted dissipation laws and conserved quantities.

The paper is organized as follows. Section 2 reviews the basic notions of $q$-contact geometry and develops the Hamiltonian formalism for uniform $q$-contact manifolds, including the Liouville-type theorem, Noether theorem, $q$-contact transformations, and the generalized Noether theorem for time-dependent systems. Section 3 extends the theory to non-uniform $q$-contact structures, introducing structure endomorphisms, channel Hamiltonian vector fields, the admissibility condition, effective Hamiltonian dynamics, the associated bracket and Noether theorem, and the evolution of volume forms. Section 4 presents the elastoplastic damage model with two internal variables as an explicit application, including the construction of the non-uniform $2$-contact structure, derivation of the effective equations of motion, and numerical validation. Finally, Section 5 concludes the paper with a summary of results and a discussion of future research directions, including constrained systems, symmetry reduction, geometric integrators, and applications to continuum mechanics and thermodynamics.
\section{Hamiltonian systems on uniform $q$-contact manifolds}

We begin by recalling some basic notions from contact geometry. 

Let $(M,\eta)$ be a contact manifold, that is, a $(2n+1)$-dimensional manifold endowed with a $1$-form $\eta$ such that
\[
\eta\wedge(d\eta)^n\neq 0.
\]
Then there exists a unique vector field $\mathcal R$, called the \emph{Reeb vector field}, satisfying
\[
i_{\mathcal R}d\eta=0,
\qquad
i_{\mathcal R}\eta=1.
\]
Equivalently, a contact structure may be described by the contact distribution
\[
\mathcal H=\ker\eta.
\]
Conversely, a contact distribution is globally defined by a contact form if and only if it is co-orientable \cite{Arnold,Le,Loose,Willett}.

A natural higher-codimensional generalization of contact geometry is provided by $q$-contact structures, which may be regarded as the geometric structures underlying contact foliations.

\begin{definition}[$q$-contact manifolds {\cite{Almeida,Finamore,Finamore2}}]\label{Defq}
Let $n,q$ be positive integers and let $M$ be a smooth manifold of dimension $2n+q$. A \emph{$q$-contact structure} on $M$ is a collection
\[
\vec{\lambda}=(\lambda_1,\ldots,\lambda_q)
\]
of pointwise linearly independent $1$-forms together with a splitting
\[
TM=\mathcal R\oplus\xi
\]
of the tangent bundle such that
\begin{enumerate}
\item
$
\xi=\bigcap_{i=1}^q\ker\lambda_i;$

\item for every $i=1,\ldots,q$, the restriction
$
d\lambda_i|_\xi
$
is non-degenerate;

\item for every $i=1,\ldots,q$,
$
\ker d\lambda_i=\mathcal R.$

\end{enumerate}
\end{definition}

\begin{remark}
Since the forms $\lambda_1,\ldots,\lambda_q$ are linearly independent, the distribution $\xi$ has constant rank $2n$. Consequently, condition (ii) is equivalent to
\[
(d\lambda_i|_\xi)^n\neq 0,
\]
or, equivalently, to $(\xi,d\lambda_i|_\xi)$ being a symplectic vector bundle over $M$.
\end{remark}

A manifold endowed with a $q$-contact structure is called a \emph{$q$-contact manifold} and is denoted by
\[
(M,\vec{\lambda},\mathcal R\oplus\xi),
\]
or simply by $M$ when no confusion can arise.

The collection $\{\lambda_i\}_{i=1}^q$ is called an \emph{adapted coframe} of the $q$-contact structure. The distributions $\mathcal R$ and $\xi$ are called the \emph{Reeb distribution} and the \emph{$q$-contact distribution}, respectively. Sections of $\xi$ are called \emph{horizontal vector fields}.

\begin{proposition}[{\cite{Almeida}}]\label{P1}
There exists a unique collection of linearly independent vector fields $R_1,\ldots,R_q$ tangent to $\mathcal R$ such that $\lambda_i(R_j)=\delta_i^j$ for $i,j=1,\ldots,q$, and $\mathcal R = \operatorname{Span}\{R_1,\ldots,R_q\}$. Moreover, $[R_i,R_j]=0$ for $i,j=1,\ldots,q$.
\end{proposition}

\begin{proof}
Let $R_1,\ldots,R_q$ be the unique vector fields satisfying $\lambda_i(R_j)=\delta_i^j$ and $\mathcal R = \operatorname{Span}\{R_1,\ldots,R_q\}$. For every $k=1,\ldots,q$, we have
\[
0 = L_{R_i}(i_{R_j}\lambda_k) = i_{[R_i,R_j]}\lambda_k + i_{R_j}(L_{R_i}\lambda_k).
\]
Since $R_i \in \ker d\lambda_k$, Cartan's formula gives $L_{R_i}\lambda_k = i_{R_i}d\lambda_k + d(i_{R_i}\lambda_k) = 0$. Therefore $i_{[R_i,R_j]}\lambda_k = 0$ for $k=1,\ldots,q$, which implies that $[R_i,R_j] \in \xi$.

On the other hand,
\[
0 = L_{R_i}(i_{R_j}d\lambda_k) = i_{[R_i,R_j]}d\lambda_k + i_{R_j}(L_{R_i}d\lambda_k).
\]
Since $L_{R_i}d\lambda_k = d(L_{R_i}\lambda_k) = 0$, it follows that $i_{[R_i,R_j]}d\lambda_k = 0$. Hence $[R_i,R_j] \in \ker d\lambda_k = \mathcal R$.

Because $[R_i,R_j] \in \xi$ and $[R_i,R_j] \in \mathcal R$, and $TM = \mathcal R \oplus \xi$, we conclude that $[R_i,R_j] = 0$.
\end{proof}
\begin{remark}
Proposition \ref{P1} implies that the Reeb distribution $\mathcal R$ is Frobenius integrable.
\end{remark}
\begin{definition}[Uniform $q$-contact structures]\label{D3}
A $q$-contact structure with adapted coframe
\[
\vec{\lambda}=(\lambda_1,\ldots,\lambda_q)
\]
is called \emph{uniform} if
\[
d\lambda_i=d\lambda_j,
\qquad
i,j=1,\ldots,q.
\]
A manifold endowed with such a structure is called a \emph{uniform $q$-contact manifold}.
\end{definition}
If $X$ is a vector field on a uniform $q$-contact manifold $M$, then it admits a unique decomposition
\[
X = \sum_{i=1}^q \lambda_i(X) R_i + \widetilde X,
\]
where $\widetilde X \in \Gamma(\xi)$ is a horizontal vector field.

Similarly, taking annihilators yields the decomposition $T^*M = \mathcal R^\circ \oplus \xi^\circ$. The elements of $\mathcal R^\circ$, namely the $1$-forms $\beta$ satisfying $\beta(R_i)=0$ for $i=1,\ldots,q$, are called \emph{semi-basic forms}. Therefore every differential $1$-form $\beta$ on $M$ can be written uniquely as
\[
\beta = \sum_{i=1}^q (i_{R_i}\beta)\,\lambda_i + \widehat\beta,
\]
where $\widehat\beta$ is a semi-basic form.

The vector fields $R_1,\ldots,R_q$ are called the \emph{Reeb vector fields} of the $q$-contact structure.

As is well known, on a contact manifold $(M,\eta)$, the Hamiltonian vector field associated with a Hamiltonian function $H \in C^\infty(M)$ is defined by
\[
\begin{cases}
\eta(X_H) = -H, \\
i_{X_H} d\eta = dH - \mathcal R(H)\eta.
\end{cases}
\]
This construction admits a natural generalization to uniform $q$-contact manifolds.

On a uniform $q$-contact manifold $(M,\vec{\lambda},\mathcal R \oplus \xi)$, every smooth function determines a unique Hamiltonian vector field (see \cite{Zhao}).

\begin{theorem}[\cite{Zhao}]\label{T2}
Let $(M,\vec{\lambda},\mathcal R \oplus \xi)$ be a uniform $q$-contact manifold and let $H \in C^\infty(M)$. Then there exists a unique vector field $X_H$ satisfying
\begin{align}\label{CH}
\begin{cases}
\lambda_i(X_H) = -H, \qquad i=1,\ldots,q, \\
i_{X_H} d\lambda_1 = dH - \displaystyle\sum_{i=1}^q dH(R_i)\,\lambda_i.
\end{cases}
\end{align}
\end{theorem}






\begin{definition}\label{DeH}
The vector field $X_H$ given by Theorem \ref{T2} is called the
\emph{$q$-contact Hamiltonian vector field} associated with $H$.

\end{definition}

\begin{remark}\label{PL}
According to the discussion in \cite{Zhao}, we know that
\eqref{CH} is equivalent to
\begin{align}\label{CH2}
\begin{cases}
\lambda_i(X_H)=-H,
\\[2mm]
L_{X_H}\lambda_i
=
-\displaystyle\sum_{j=1}^qdH(R_j)\lambda_j,
\qquad i=1,\ldots,q.
\end{cases}
\end{align}
\end{remark}

We now proceed to define the evolution vector field associated with the $q$-contact Hamiltonian vector field.
\begin{definition}\label{DeE}
The vector field $E_H$ defined by
\begin{align}\label{EH}
E_H := X_H + \sum_{i=1}^q H \, R_i
\end{align}
is called the \emph{evolution vector field} associated with $H$, where $X_H$ is the $q$-contact Hamiltonian vector field given by Theorem \ref{T2} and $R_i$ are the Reeb vector fields corresponding to the $q$-contact structure $(\vec{\lambda},\mathcal R \oplus \xi)$.
\end{definition}
\begin{remark}
Unlike the $q$-contact Hamiltonian vector field $X_H$, the evolution vector field $E_H$ belongs to the kernel of the $q$-contact forms, i.e.,
\begin{align}
\lambda_i(E_H) = 0, \qquad i=1,\ldots,q,
\end{align}
since $\lambda_i(X_H) = -H$ and $\lambda_i(R_j) = \delta_{ij}$. In the context of thermodynamics, this property makes $E_H$ particularly relevant, as its integral curves naturally describe the evolution of equilibrium states, which are represented by Legendre submanifolds. Thus, $E_H$ can be interpreted as a deformation of $X_H$ along the Reeb directions.
\end{remark}

We now give the bracket associated with the $q$-contact Hamiltonian dynamics, which is introduced in \cite{Zhao}.
\begin{definition}\label{D2}
Let $(M,\vec{\lambda},\mathcal R \oplus \xi)$ be a uniform $q$-contact manifold. The \emph{$q$-contact bracket} is the map
\[
\{\cdot,\cdot\}: C^\infty(M) \times C^\infty(M) \longrightarrow C^\infty(M)
\]
defined by
\begin{align}\label{dissipation}
\{f,g\}
= \lambda_1([X_f,X_g]) = L_{X_f}(\lambda_1(X_g)) - (L_{X_f}\lambda_1)(X_g) = -X_f(g) - g\sum_{i=1}^q R_i(f). 
\end{align}
\end{definition}

\begin{remark}\label{R3}
The definition of the bracket is independent of the choice of the contact form. Indeed, since
\[
d\lambda_i(X,Y)
=
X(\lambda_i(Y))
-
Y(\lambda_i(X))
-
\lambda_i([X,Y]),
\]
and $\lambda_i(X_f)=\lambda_i(X_g)=0$, we obtain
\[
\lambda_i([X_f,X_g])
=
-d\lambda_i(X_f,X_g).
\]
Because the structure is uniform,
\[
d\lambda_1=\cdots=d\lambda_q,
\]
it follows that
\[
\lambda_i([X_f,X_g])
=
\lambda_j([X_f,X_g]),
\qquad
i,j=1,\ldots,q.
\]
\end{remark}
It is well known that every uniform $q$-contact manifold admits local coordinates
\begin{equation}\label{qcontact}
(x_1,\ldots,x_n,y_1,\ldots,y_n,z_1,\ldots,z_q),
\end{equation}
called \emph{$q$-contact coordinates}, in which the adapted coframe and the Reeb vector fields are given by \cite{Blair}
\[
\lambda_i
=
dz_i-\sum_{a=1}^n x_a\,dy_a,
\qquad
R_i
=
\frac{\partial}{\partial z_i},
\qquad
i=1,\ldots,q.
\]

In these coordinates, the Hamiltonian vector field associated with a Hamiltonian function $H$ takes the form
\begin{equation}\label{XH}
X_H
=
\sum_{i=1}^{q}
\left(
x_a\frac{\partial H}{\partial x_a}
-
H
\right)
\frac{\partial}{\partial z_i}
-
\sum_{a=1}^{n}
\left(
x_a\sum_{i=1}^{q}\frac{\partial H}{\partial z_i}
+
\frac{\partial H}{\partial y_a}
\right)
\frac{\partial}{\partial x_a}
+
\sum_{a=1}^{n}
\frac{\partial H}{\partial x_a}
\frac{\partial}{\partial y_a}.
\end{equation}

Therefore, the corresponding Hamiltonian equations are
\begin{align}
\dot y_a
&=
\frac{\partial H}{\partial x_a},
\label{H11}
\\
\dot x_a
&=
-
x_a\sum_{i=1}^{q}\frac{\partial H}{\partial z_i}
-
\frac{\partial H}{\partial y_a},
\label{H22}
\\
\dot z_i
&=
x_a\frac{\partial H}{\partial x_a}
-
H.
\label{H3}
\end{align}

\begin{remark}
Equations \eqref{H11}--\eqref{H3} generalize both the Hamiltonian equations on symplectic manifolds and the Hamiltonian equations on contact manifolds.
Indeed, when $q=1$, they reduce to the standard contact Hamiltonian equations.

Furthermore, if the Hamiltonian is independent of the variables $z_1,\ldots,z_q$, then
\[
\frac{\partial H}{\partial z_i}=0,
\qquad
i=1,\ldots,q,
\]
and equations \eqref{H11}--\eqref{H22} reduce to the classical Hamilton equations on a symplectic manifold. In this case, the Hamiltonian function is conserved along the flow.

A more general, although non-generic, situation occurs when
\[
\sum_{i=1}^{q}\frac{\partial H}{\partial z_i}=0.
\]
In this case, the dissipative contribution in equation \eqref{H22} vanishes due to an exact cancellation between the different dissipation channels, and the $(x_a,y_a)$-dynamics is again governed by the classical symplectic Hamiltonian equations.
\end{remark}
\begin{example}
Consider the Hamiltonian fucntion
\[
H
=
\frac{x^2}{2m}
+
V(y)
+
\sum_{i=1}^{q}\gamma_i z_i,
\]
where $V(y)$ is a potential function and
\[
\gamma_i\in\mathbb R,
\qquad
i=1,\ldots,q,
\]
are constants.

The Hamiltonian equations become
\begin{align}
\dot y
&=
\frac{x}{m},
\label{E1}
\\
\dot x
&=
-\frac{\partial V}{\partial y}
-
x(\gamma_1+\cdots+\gamma_q),
\label{E2}
\\
\dot z_i
&=
\frac{x^2}{2m}
-
V(y)
-
\sum_{j=1}^{q}\gamma_j z_j.
\label{E3}
\end{align}

Combining \eqref{E1} and \eqref{E2}, we obtain
\[
\ddot y
+
(\gamma_1+\cdots+\gamma_q)\dot y
+
\frac1m\frac{\partial V}{\partial y}
=
0,
\]
which is the equation of motion of a damped mechanical system with damping coefficient
\[
\gamma_1+\cdots+\gamma_q.
\]
\end{example}
\subsection{Dynamics on uniform $q$-contact manifolds}
We now investigate the evolution of observables along the Hamiltonian flow.

Let $F \in C^\infty(M)$. Using Definition \ref{D2}, we obtain
\[
\{H,F\} = -X_H(F) - F\sum_{i=1}^{q} R_i(H).
\]
Equivalently,
\begin{equation}\label{evolution}
X_H(F) = -\{H,F\} - F\sum_{i=1}^{q} R_i(H).
\end{equation}
Since $\frac{dF}{dt} = X_H(F)$, equation \eqref{evolution} describes the evolution of any observable along the Hamiltonian flow.

\begin{definition}
Let $X_H$ be the Hamiltonian vector field associated with a Hamiltonian function $H$.
\begin{enumerate}
\item A smooth function $f \in C^\infty(M)$ is called a \emph{dissipated quantity} if $X_H(f) = -f\sum_{i=1}^q R_i(H)$.
\item A smooth function $g \in C^\infty(M)$ is called a \emph{constant of motion} (or \emph{conserved quantity}) if $X_H(g) = 0$.
\end{enumerate}
\end{definition}

\begin{remark}
Setting $F=H$ in the evolution equation $X_H(F) = -\{H,F\} - F\sum_{i=1}^q R_i(H)$, and using the skew-symmetry of the $q$-contact bracket, we obtain
\[
\frac{dH}{dt} = X_H(H) = -H\sum_{i=1}^q R_i(H).
\]
Hence the Hamiltonian function is always a dissipated quantity. Moreover, $H$ is a constant of motion if and only if $H\sum_{i=1}^q R_i(H) = 0$.

In particular, in $q$-contact coordinates, $R_i = \partial_{z_i}$, and therefore $\sum_{i=1}^q R_i(H) = \sum_{i=1}^q \partial H/\partial z_i$. Consequently, a sufficient condition for $H$ to be conserved is $\sum_{i=1}^q \partial H/\partial z_i = 0$. When $H$ is independent of the variables $(z_1,\ldots,z_q)$, this condition is automatically satisfied.
\end{remark}
Now that we have introduced the evolution vector field $E_H$ in the previous definition, we proceed to state a proposition which establishes that the dissipative and conservative quantities of the Hamiltonian vector field $X_H$ admit an alternative characterization via the evolution vector field $E_H$.
\begin{proposition}
A smooth function $f$ is a dissipated quantity if and only if
\[
E_H(f) = \sum_{i=1}^q \big( H R_i(f) - f R_i(H) \big).
\]
In particular, if $f$ is nonvanishing, then $f$ is a dissipated quantity if and only if
\[
E_H\!\left(\frac{1}{f}\right) = \sum_{i=1}^q R_i\!\left(\frac{H}{f}\right).
\]
Moreover, $f$ is a conservative quantity if and only if
\[
E_H(f) = \sum_{i=1}^q H R_i(f).
\]
\end{proposition}
\begin{proof}
The proof is straightforward.
\end{proof}
\subsection{Uniform $q$-contact transformations}
We now introduce the analogue of canonical transformations in the uniform $q$-contact setting.

\begin{definition}
A diffeomorphism $\Phi: M \to M$ is called a \emph{$q$-contact transformation} if there exists a non-zero constant $f \in \mathbb R$ such that $\Phi^*\lambda_i = f\lambda_i$ for $i=1,\ldots,q$.
\end{definition}

\begin{remark}
Since $d\lambda_1 = \cdots = d\lambda_q$, it follows that $\Phi^*(d\lambda_i) = f\,d\lambda_1$ for $i=1,\ldots,q$.
\end{remark}

In $q$-contact coordinates, a coordinate transformation
\[
(x_a, y_a, z_i) \longmapsto (\widetilde x_a, \widetilde y_a, \widetilde z_i)
\]
is a $q$-contact transformation if and only if
\begin{equation}\label{fl}
f(dz_i - x_a\,dy_a) = d\widetilde z_i - \widetilde x_a\,d\widetilde y_a, \qquad i=1,\ldots,q.
\end{equation}

Expanding \eqref{fl} yields
\begin{align*}
f &= \frac{\partial\widetilde z_1}{\partial z_1} - \sum_{a=1}^n \widetilde x_a \frac{\partial\widetilde y_a}{\partial z_1} = \cdots = \frac{\partial\widetilde z_q}{\partial z_q} - \sum_{a=1}^n \widetilde x_a \frac{\partial\widetilde y_a}{\partial z_q}, \\
-fx_j &= \frac{\partial\widetilde z_i}{\partial y_j} - \sum_{a=1}^n \widetilde x_a \frac{\partial\widetilde y_a}{\partial y_j}, \\
0 &= \frac{\partial\widetilde z_i}{\partial x_j} - \sum_{a=1}^n \widetilde x_a \frac{\partial\widetilde y_a}{\partial x_j}, \\
0 &= \frac{\partial\widetilde z_i}{\partial z_j}, \qquad i \neq j.
\end{align*}

Assuming that $(y_a, \widetilde y_a, z_i)$ are independent coordinates, we may introduce generating functions $\widetilde z_i = \widetilde z_i(y_a, \widetilde y_a, z_i)$. Their differentials satisfy
\begin{equation}\label{TS}
d\widetilde z_i = \frac{\partial\widetilde z_i}{\partial z_i}\,dz_i + \sum_{a=1}^n \frac{\partial\widetilde z_i}{\partial y_a}\,dy_a + \sum_{a=1}^n \frac{\partial\widetilde z_i}{\partial \widetilde y_a}\,d\widetilde y_a.
\end{equation}

Substituting \eqref{TS} into \eqref{fl} yields
\begin{equation}\label{ZH}
f = \frac{\partial\widetilde z_i}{\partial z_i}, \qquad fx_j = -\frac{\partial\widetilde z_i}{\partial y_j}, \qquad \widetilde x_a = \frac{\partial\widetilde z_i}{\partial \widetilde y_a}.
\end{equation}

In the particular case $f=1$, one has
\[
\widetilde z_i=z_i-F(y_a,\widetilde y_a),
\]
where $F$ generates the induced transformation on the horizontal symplectic distribution. The corresponding relations
\[
x_j=\frac{\partial F}{\partial y_j},
\qquad
\widetilde x_j=-\frac{\partial F}{\partial\widetilde y_j},
\]
are formally identical to those of a generating function for a symplectic canonical transformation, although here they arise from the symplectic structure on the horizontal distribution rather than from the canonical symplectic structure of a cotangent bundle.
\subsection{Liouville-type theorem and symmetries on uniform $q$-contact manifolds}
Using our previous discussion, we can readily observe the following result.
\begin{proposition}[Liouville-type formula]
Let $\mu = \lambda_1 \wedge \cdots \wedge \lambda_q \wedge (d\lambda_1)^n$. Then the Hamiltonian flow satisfies
\begin{equation}\label{LX}
L_{X_H}\mu = -(n+1)\left( \sum_{i=1}^q R_i(H) \right)\mu.
\end{equation}
\end{proposition}

\begin{corollary}
The divergence of the Hamiltonian vector field with respect to the volume form $\mu$ is
\begin{equation}\label{LXX}
\operatorname{div}_{\mu}(X_H) = -(n+1)\sum_{i=1}^q R_i(H).
\end{equation}
\end{corollary}
By using the above results, we obtain the following theorem.
\begin{theorem}[Liouville-type theorem]
Let $X_H$ be the Hamiltonian vector field associated with a Hamiltonian function $H$ on a uniform $q$-contact manifold $(M,\vec{\lambda},\mathcal R \oplus \xi)$. Then $X_H$ is divergence-free with respect to the volume form $\mu = \lambda_1 \wedge \cdots \wedge \lambda_q \wedge (d\lambda_1)^n$ if and only if $\sum_{i=1}^q R_i(H) = 0$.

Moreover, if $H$ is nowhere vanishing, then $\mu_H := H^{-(n+1)}$ is a Jacobi multiplier of $X_H$, that is, $\mu_H X_H$ is divergence-free with respect to the volume form $\mu$.
\end{theorem}

\begin{proof}
The first statement follows directly from equation \eqref{LXX}.

For the second statement, let $\mu_H = H^{-(n+1)}$. Using $X_H(H) = -H\sum_{i=1}^q R_i(H)$ together with \eqref{LX}, we obtain
\begin{align*}
L_{\mu_H X_H}\mu
= L_{X_H}(\mu_H \mu) = (L_{X_H}\mu_H)\mu + \mu_H L_{X_H}\mu.
\end{align*}
Since $L_{X_H}\mu_H = -(n+1)H^{-(n+2)}X_H(H)$, it follows that $L_{X_H}\mu_H = (n+1)\mu_H \sum_{i=1}^q R_i(H)$. Combining this with \eqref{LX}, we obtain
\[
L_{\mu_H X_H}\mu
= (n+1)\mu_H\left( \sum_{i=1}^q R_i(H) \right)\mu - (n+1)\mu_H\left( \sum_{i=1}^q R_i(H) \right)\mu = 0.
\]
Therefore $\mu_H X_H$ is divergence-free with respect to $\mu$.
\end{proof}

\begin{remark}
If $R_i(H) = 0$ for $i=1,\ldots,q$, then $L_{X_H}\lambda_i = 0$ for $i=1,\ldots,q$, and therefore the flow of $X_H$ consists of $q$-contact transformations preserving the uniform $q$-contact structure. More generally, the weaker condition $\sum_{i=1}^q R_i(H) = 0$ already guarantees preservation of the volume form $\lambda_1 \wedge \cdots \wedge \lambda_q \wedge (d\lambda_1)^n$.
\end{remark}

\begin{definition}\label{Noth}
A Hamiltonian vector field $X_F$ is called a \emph{Noether symmetry} of the Hamiltonian vector field $X_H$ if
\[
\{F,H\}=0.
\]
\end{definition}
The corresponding $q$-contact version of Noether's theorem is the following.
\begin{theorem}[$q$-contact Noether theorem]
A Hamiltonian vector field $X_F$ is a Noether symmetry of $X_H$ if and only if its Hamiltonian function $F = -i_{X_F}\lambda_1$ is a dissipated quantity.
\end{theorem}

\begin{proof}
Using $F = -i_{X_F}\lambda_1$, we obtain
\[
L_{X_H}F = -L_{X_H}(i_{X_F}\lambda_1) = -i_{[X_H,X_F]}\lambda_1 - i_{X_F}(L_{X_H}\lambda_1).
\]
Since $X_F$ is a Noether symmetry, $\{F,H\} = \lambda_1([X_F,X_H]) = 0$, and therefore $i_{[X_H,X_F]}\lambda_1 = 0$.

Using Remark \ref{PL}, $L_{X_H}\lambda_1 = -\sum_{j=1}^q R_j(H)\,\lambda_j$. Hence
\[
L_{X_H}F = i_{X_F}\left( \sum_{j=1}^q R_j(H)\,\lambda_j \right) = -F\sum_{j=1}^q R_j(H),
\]
which is precisely the definition of a dissipated quantity.
\end{proof}
\begin{definition}
Let $X_H$ be a Hamiltonian vector field on a uniform $q$-contact manifold. A vector field $Y \in \mathfrak X(M)$ is called a \emph{dynamical symmetry} of $X_H$ if there exists a smooth function $\Lambda \in C^\infty(M)$ such that $[Y,X_H] = \Lambda X_H$. When $\Lambda = 0$, the symmetry is called \emph{trivial}.
\end{definition}

\begin{theorem}
Let $Y = X_F$ be a Hamiltonian dynamical symmetry of $X_H$. If $[Y,X_H] = \Lambda X_H$ with $\Lambda \neq 0$, then $Y$ is not a Noether symmetry.
\end{theorem}

\begin{proof}
Since $[Y,X_H] = \Lambda X_H$, we obtain
\[
i_{[Y,X_H]}\lambda_1 = \Lambda\, i_{X_H}\lambda_1 = -\Lambda H.
\]
Assume by contradiction that $Y = X_F$ is a Noether symmetry. Then $\{F,H\} = \lambda_1([X_F,X_H]) = \lambda_1([Y,X_H]) = 0$. However,
\[
\lambda_1([Y,X_H]) = i_{[Y,X_H]}\lambda_1 = -\Lambda H.
\]
Since $\Lambda \neq 0$, this contradicts the previous equality whenever $H \not\equiv 0$. Therefore $Y$ cannot be a Noether symmetry.
\end{proof}
We now consider an extended \(q\)-contact phase space as the manifold \(M \times \mathbb{R}\), endowed with the \(1\)-forms
\[
\lambda_i^E = \pi_t^* \lambda_i + H\,dt,
\]
where \(t\) is the coordinate on \(\mathbb{R}\), and \(\pi_t : M \times \mathbb{R} \to \mathbb{R}\) denotes the canonical projection. For simplicity, we shall write this directly as
\[
\lambda_i^E = \lambda_i + H\,dt.
\]
Assume that $R_i(H) = 0$ for $i=1,\ldots,q$. Then the defining equations
\[
i_{X_H^t} d\lambda_i^E = 0, \qquad i_{X_H^t}\lambda_i^E = 0,
\]
do not determine $X_H^t$ uniquely. Indeed, if $X_H^t$ is a solution, then so is $f X_H^t$ for any nowhere-vanishing function $f: M \times \mathbb R \to \mathbb R \setminus \{0\}$. Therefore, $X_H^t$ is defined only up to a nowhere-vanishing rescaling.

Choosing the normalization $f=1$, we recover in local coordinates $$(x_1,\ldots,x_n,y_1,\ldots,y_n,z_1,\ldots,z_q,t)$$ the time-dependent $q$-contact Hamiltonian equations
\begin{align}
\dot y_a &= \frac{\partial H}{\partial x_a}, \\
\dot x_a &= -x_a\sum_{i=1}^q \frac{\partial H}{\partial z_i} - \frac{\partial H}{\partial y_a}, \\
\dot z_i &= x_a\frac{\partial H}{\partial x_a} - H, \\
\dot t &= 1.
\end{align}
Hence changing the function $f$ only amounts to a reparametrization of the trajectories.

\begin{remark}
The construction presented above bears some resemblance to the notion of time-dependent contact manifolds introduced in \cite{deLeonGasetGraciaMunozRivas2023}, where contact geometry is extended by incorporating an explicit time variable. Although our setting is fundamentally different, involving multiple contact forms rather than an additional time coordinate, it would be interesting to investigate whether the present framework admits a unified formulation in terms of a suitable notion of time-dependent $q$-contact manifolds. We leave this question for future work.
\end{remark}

To simplify the notation, we shall assume from now on that $f=1$. Thus $X_H^t = X_H + \partial_t$, where $X_H$ denotes the $q$-contact Hamiltonian vector field associated with $H$.

We now investigate symmetries and dissipated quantities in the extended phase space.

\begin{definition}
A vector field $Y \in \mathfrak X(M \times \mathbb R)$ is called a \emph{generalized Noether symmetry} if $L_Y \lambda_i^E = \sigma_i \lambda_i^E$ for $i=1,\ldots,q$, for some functions $\sigma_i \in C^\infty(M \times \mathbb R)$.
\end{definition}
\begin{proposition}
Generalized Noether symmetries form a Lie algebra with respect to the Lie bracket of vector fields.
\end{proposition}

\begin{proof}
Let $L_{Y_1}\lambda_i^E = \sigma_i^1\lambda_i^E$ and $L_{Y_2}\lambda_i^E = \sigma_i^2\lambda_i^E$. Then
\[
L_{[Y_1,Y_2]}\lambda_i^E
= L_{Y_1}L_{Y_2}\lambda_i^E - L_{Y_2}L_{Y_1}\lambda_i^E
= \bigl( Y_1(\sigma_i^2) - Y_2(\sigma_i^1) \bigr)\lambda_i^E.
\]
Hence $[Y_1,Y_2]$ is again a generalized Noether symmetry.
\end{proof}
\begin{definition}
A smooth function $F: M \times \mathbb R \to \mathbb R$ is called a \emph{dissipated quantity} in the extended phase space if
\begin{equation}\label{Diss}
L_{X_H^t} F = -F \sum_{i=1}^q R_i(H).
\end{equation}
\end{definition}
Equation \eqref{Diss} will be referred to as the \emph{dissipation equation}.
Observe that
\[
L_{X_H^t} H = -H \sum_{i=1}^q R_i(H) + \frac{\partial H}{\partial t}.
\]
Consequently, when $H$ depends explicitly on time, the Hamiltonian function is generally not a dissipated quantity.

Moreover, in the time-independent theory, if $F$ is dissipated and $H$ is nowhere vanishing, then $F/H$ is conserved. This property no longer holds in general when $H$ depends on time. Indeed,
\[
L_{X_H^t}\left( \frac{F}{H} \right) = -\frac{F}{H^2} \frac{\partial H}{\partial t}.
\]
Therefore $F/H$ is conserved precisely when $\partial H/\partial t = 0$.

We firstly establish a preliminary result.
We are now ready to prove the generalized Noether theorem.

\begin{theorem}[Generalized Noether theorem]
Let $Y$ be a generalized Noether symmetry of a time-dependent $q$-contact Hamiltonian system. Then the functions $F_i = i_Y\lambda_i^E$ for $i=1,\ldots,q$ are dissipated quantities.
\end{theorem}

\begin{proof}
For every $i=1,\ldots,q$,
\begin{align*}
L_{X_H^t}F_i
&= L_{X_H^t}(i_Y\lambda_i^E) \\
&= i_{[X_H^t,Y]}\lambda_i^E + i_Y(L_{X_H^t}\lambda_i^E).
\end{align*}
The previous proposition implies $i_{[X_H^t,Y]}\lambda_i^E = 0$. Moreover,
\[
L_{X_H^t}\lambda_i^E = -\left( \sum_{j=1}^q R_j(H) \right)\lambda_i^E.
\]
Hence
\[
L_{X_H^t}F_i = -\left( \sum_{j=1}^q R_j(H) \right) i_Y\lambda_i^E = -\left( \sum_{j=1}^q R_j(H) \right) F_i.
\]
Therefore $F_i$ satisfies the dissipation equation.
\end{proof}
\section{Hamiltonian systems on non-uniform $q$-contact structures}

In many physical and geometric applications, the uniformity condition
\[
d\lambda_i = d\lambda_1, \qquad i = 1, \ldots, q,
\]
is overly restrictive. Distinct dissipative mechanisms may interact with the phase space through different geometric channels, naturally leading to situations in which the exterior derivatives of the adapted coframe are no longer identical.

We now return to the most general definition of a $q$-contact manifold, namely Definition~\ref{Defq}, which does not impose the condition $d\lambda_i = d\lambda_j$ required in the case of a uniform $q$-contact manifold. To distinguish this more general setting from the uniform case, we refer to such a $q$-contact manifold as a non-uniform $q$-contact manifold.

For the reader's convenience, we restate the content of Definition~\ref{Defq}, where the manifold is now referred to as a non-uniform $q$-contact manifold.
\begin{definition}[non-uniform $q$-contact manifolds]
Let $n,q$ be positive integers and let $M$ be a smooth manifold of dimension $2n+q$. A \emph{non-uniform $q$-contact structure} on $M$ is a collection
\[
\vec{\lambda}=(\lambda_1,\ldots,\lambda_q)
\]
of pointwise linearly independent $1$-forms together with a splitting
\[
TM=\mathcal R\oplus\xi
\]
of the tangent bundle such that
\begin{enumerate}
\item
$
\xi=\bigcap_{i=1}^q\ker\lambda_i;$

\item for every $i=1,\ldots,q$, the restriction
$
d\lambda_i|_\xi
$
is non-degenerate;

\item for every $i=1,\ldots,q$,
$
\ker d\lambda_i=\mathcal R.$

\end{enumerate}
\end{definition}




\begin{remark}
Each restriction
\[
d\lambda_i|_\xi
\]
defines a symplectic structure on the horizontal distribution $\xi$.
Unlike the uniform case, these symplectic structures need not coincide.
Thus the horizontal geometry is governed by a family of symplectic forms rather than by a single one.
\end{remark}





We now investigate the relations among $d\lambda_i$, $i = 1, \dots, q$.
Fix the symplectic form
\[
\omega:=d\lambda_1|_\xi.
\]
We can give the following proposition.
\begin{proposition}
For each $i = 1,\dots,q$, there exists a unique vector bundle endomorphism
\[
B_i: \xi \to \xi
\]
such that for all $X, Y \in \xi$,
\[
d\lambda_i(X,Y) = \omega(B_i X, Y).
\]
\end{proposition}

\begin{proof}
Since $\omega$ is non-degenerate on $\xi$, the bundle morphism
\[
\omega^\flat: \xi \to \xi^*, \quad X \mapsto \omega(X, \cdot)
\]
is an isomorphism, where $\xi^*$ denotes the dual bundle of $\xi$.

Fix $X \in \xi$. Then the map $Y \mapsto d\lambda_i(X,Y)$ defines an element of $\xi^*$. Hence there exists a unique vector $B_i X \in \xi$ such that
\[
\omega(B_i X, Y) = d\lambda_i(X,Y)
\]
for every $Y \in \xi$.

The assignment $X \mapsto B_i X$ is smooth and linear, and therefore defines a unique vector bundle endomorphism $B_i: \xi \to \xi$.
\end{proof}

\begin{remark}
The endomorphisms $B_i = \omega^{-1} \circ d\lambda_i|_\xi$ measure the deviation of the symplectic structures $d\lambda_i|_\xi$ from the reference symplectic form $\omega = d\lambda_1|_\xi$. In particular, $B_1 = \operatorname{Id}_\xi$.
\end{remark}

\begin{remark}
The above construction depends on the choice of the reference symplectic form
\(
\omega=d\lambda_1|_\xi.
\)
More intrinsically, given any pair of indices $i,j\in\{1,\ldots,q\}$, one may define the vector bundle automorphism
\[
B_{ij}
=
(d\lambda_j^\flat)^{-1}\circ d\lambda_i^\flat,
\]
where
\[
d\lambda_i^\flat(X)=d\lambda_i(X,\cdot).
\]
These satisfy
\[
d\lambda_i(X,Y)
=
d\lambda_j(B_{ij}X,Y),
\]
for all $X,Y\in\xi$. The family $\{B_{ij}\}$ provides a coordinate-free comparison among the different symplectic structures induced on $\xi$. In the present work, however, fixing a reference form considerably simplifies the exposition and is sufficient for the development of the subsequent theory.
\end{remark}

\begin{definition}
The bundle morphisms $B_i: \xi \to \xi$ are called the \emph{structure endomorphisms} associated with the non-uniform $q$-contact structure.
\end{definition}
According to Proposition~\ref{P1}, there exists a unique collection of linearly independent vector fields $R_1,\ldots,R_q$ tangent to $\mathcal R$ such that $\lambda_i(R_j)=\delta_i^j$ for $i,j=1,\ldots,q$, and $\mathcal R = \operatorname{Span}\{R_1,\ldots,R_q\}$. Moreover, $[R_i,R_j]=0$ for $i,j=1,\ldots,q$. Using these $R_i$, we now define the Hamiltonian dynamics associated with the different geometric channels.
\begin{theorem}\label{THchannel}
Let $(M,\vec\lambda,\mathcal R\oplus\xi)$ be a non-uniform $q$-contact manifold and let $H\in C^\infty(M)$. For each $i=1,\ldots,q$, there exists a unique vector field $X_H^i$ satisfying
\[
\lambda_j(X_H^i) = -H, \qquad j=1,\ldots,q,
\]
and
\[
i_{X_H^i} d\lambda_i = dH - \sum_{j=1}^q dH(R_j)\,\lambda_j.
\]
\end{theorem}

\begin{proof}
Define the $1$-form $\beta = dH - \sum_{j=1}^q dH(R_j)\,\lambda_j$. For every $k=1,\ldots,q$,
\[
\beta(R_k) = dH(R_k) - \sum_{j=1}^q dH(R_j)\,\lambda_j(R_k).
\]
Using $\lambda_j(R_k) = \delta_j^k$, we obtain $\beta(R_k)=0$. Hence $\beta$ annihilates $\mathcal R$ and therefore may be regarded as a section of $\xi^*$.

Since $d\lambda_i|_\xi$ is non-degenerate, there exists a unique horizontal vector field $Y_H^i \in \Gamma(\xi)$ such that $i_{Y_H^i} d\lambda_i = \beta$.

Define $X_H^i = Y_H^i - H\sum_{j=1}^q R_j$. Since $Y_H^i \in \xi$, we have $\lambda_j(Y_H^i)=0$ for $j=1,\ldots,q$. Therefore
\[
\lambda_j(X_H^i) = -H\,\lambda_j\!\left( \sum_{k=1}^q R_k \right) = -H,
\qquad j=1,\ldots,q.
\]

Moreover,
\[
i_{X_H^i} d\lambda_i = i_{Y_H^i} d\lambda_i - H\sum_{j=1}^q i_{R_j} d\lambda_i.
\]
Since $R_j \in \ker d\lambda_i$, it follows that $i_{X_H^i} d\lambda_i = \beta = dH - \sum_{j=1}^q dH(R_j)\,\lambda_j$. Thus $X_H^i$ satisfies the required conditions.

Uniqueness follows from the decomposition $TM = \xi \oplus \mathcal R$ and the non-degeneracy of $d\lambda_i|_\xi$.
\end{proof}

\begin{definition}
For each $i=1,\ldots,q$, the vector field $X_H^i$ given by Theorem \ref{THchannel} is called the \emph{$i$-th channel Hamiltonian vector field} associated with the Hamiltonian function $H$.
\end{definition}

\begin{remark}
The vector fields $X_H^1,\ldots,X_H^q$ share the same Reeb component $-H\sum_{j=1}^q R_j$, and differ only through their horizontal components $Y_H^i$. Consequently, the non-uniform geometry naturally gives rise to a family of Hamiltonian dynamics, one for each symplectic structure $d\lambda_i|_\xi$.
\end{remark}
\begin{remark}
The family of channel Hamiltonian vector fields
\[
\{X_H^1,\ldots,X_H^q\}
\]
may be regarded as a $q$-vector field in the sense that the non-uniform $q$-contact structure naturally associates one Hamiltonian vector field to each symplectic structure $d\lambda_i|_\xi$. This viewpoint bears some resemblance to the $k$-symplectic formalism for classical field theories, where one considers a family of symplectic structures on the Whitney sum of $k$ copies of a cotangent bundle (see \cite{deLeonSalgadoVilarino}). Although our geometric setting is different, the coexistence of multiple symplectic structures suggests possible connections between non-uniform $q$-contact geometry and the geometric formulations of classical field theory. We leave a detailed investigation of this relationship for future work.
\end{remark}

\begin{definition}
The vector-valued field
\[
\mathbb X_H = (X_H^1,\dots,X_H^q)
\]
is called the \emph{non-uniform Hamiltonian field} associated with the Hamiltonian function $H$.
\end{definition}

\begin{remark}
The vector-valued field $\mathbb X_H = (X_H^1,\dots,X_H^q)$ describes the Hamiltonian response of the system through the different geometric channels $d\lambda_1,\dots,d\lambda_q$. Thus a single Hamiltonian function gives rise to several dynamical evolutions determined by the underlying non-uniform geometry.

Using Cartan's identity $L_X\alpha = i_X d\alpha + d(i_X\alpha)$, we obtain
\begin{align*}
L_{X_H^i}\lambda_i
&= i_{X_H^i}d\lambda_i + d(i_{X_H^i}\lambda_i) \\
&= dH - \sum_{j=1}^q dH(R_j)\,\lambda_j - dH \\
&= -\sum_{j=1}^q dH(R_j)\,\lambda_j.
\end{align*}

For $j \neq i$, we similarly obtain $L_{X_H^i}\lambda_j = i_{X_H^i}d\lambda_j - dH$. Since the $2$-forms $d\lambda_j$ and $d\lambda_i$ are generally different, no further simplification is available in general.

In the uniform case, $d\lambda_i = d\lambda_1$ for $i=1,\dots,q$, all components of the non-uniform Hamiltonian field are governed by the same symplectic structure. Consequently, the Hamiltonian vector fields $X_H^1,\dots,X_H^q$ coincide, and the non-uniform Hamiltonian dynamics reduces to the uniform $q$-contact Hamiltonian dynamics.
\end{remark}
Although the non-uniform structure naturally produces a vector-valued Hamiltonian field, it is sometimes possible to combine all geometric channels into a single effective dynamics.
\begin{definition}
Define the $2$-form $\Omega = \frac1q \sum_{i=1}^q d\lambda_i$. The non-uniform $q$-contact structure is called \emph{admissible} if $\Omega|_\xi$ is non-degenerate.
\end{definition}

Under the admissibility condition, the different geometric channels can be combined into a single effective Hamiltonian dynamics.

\begin{theorem}\label{TT}
Let $(M,\vec\lambda,\mathcal R\oplus\xi)$ be an admissible non-uniform $q$-contact manifold. Then for every $H \in C^\infty(M)$, there exists a unique vector field $X_H$ satisfying
\begin{align}
\lambda_i(X_H) &= -H, \qquad i=1,\dots,q, \label{H1} \\
i_{X_H}\Omega &= dH - \sum_{i=1}^q dH(R_i)\,\lambda_i. \label{H2}
\end{align}
\end{theorem}

\begin{proof}
Define $\beta = dH - \sum_{i=1}^q dH(R_i)\,\lambda_i$. For every $j=1,\dots,q$,
\[
\beta(R_j) = dH(R_j) - \sum_{i=1}^q dH(R_i)\,\lambda_i(R_j) = 0.
\]
Hence $\beta|_{\mathcal R}=0$, and therefore $\beta \in \xi^*$. Since $\Omega|_\xi$ is non-degenerate, there exists a unique vector field $Y_H \in \Gamma(\xi)$ such that $i_{Y_H}\Omega = \beta$.

Define $X_H = Y_H - H\sum_{i=1}^q R_i$. Since $Y_H \in \xi$, we have $\lambda_i(Y_H)=0$ for $i=1,\ldots,q$. Therefore $\lambda_i(X_H) = -H$ for $i=1,\ldots,q$.

Moreover,
\[
i_{X_H}\Omega = i_{Y_H}\Omega - H\sum_{i=1}^q i_{R_i}\Omega.
\]
Because $R_i \in \ker d\lambda_j$ for $i,j=1,\ldots,q$, it follows that $i_{R_i}\Omega = 0$. Hence
\[
i_{X_H}\Omega = \beta = dH - \sum_{i=1}^q dH(R_i)\,\lambda_i.
\]
Thus $X_H$ satisfies \eqref{H1}--\eqref{H2}.

Uniqueness follows from the decomposition $TM = \xi \oplus \mathcal R$ and the non-degeneracy of $\Omega|_\xi$.
\end{proof}
\begin{remark}
In the uniform case,
\[
d\lambda_1=\cdots=d\lambda_q,
\]
we have
\[
\Omega=d\lambda_1.
\]

Consequently, the Hamiltonian vector field provided by Theorem \ref{TT} coincides with the Hamiltonian vector field introduced in the uniform theory. Therefore, the admissible non-uniform formalism extends the uniform $q$-contact Hamiltonian framework.
\end{remark}

\begin{remark}
The vector field $X_H$ may be interpreted as the effective Hamiltonian dynamics obtained by combining all geometric channels into the averaged structure $\Omega = \frac{1}{q}\sum_{i=1}^q d\lambda_i$.

Furthermore, every uniform $q$-contact structure is admissible. Indeed, if $d\lambda_1 = \cdots = d\lambda_q$, then
\[
\Omega = \frac{1}{q}\sum_{i=1}^q d\lambda_i = d\lambda_1.
\]
Consequently, the Hamiltonian vector field defined in Theorem~\ref{TT} coincides with the Hamiltonian vector field of the uniform $q$-contact theory. Therefore, the admissible non-uniform framework naturally extends the uniform one.
\end{remark}

\begin{remark}
For the effective Hamiltonian vector field $X_H$, Cartan's identity yields
\[
L_{X_H}\Omega
=
d(i_{X_H}\Omega)
+
i_{X_H}(d\Omega).
\]

Since $\Omega$ is closed, we obtain
\[
L_{X_H}\Omega
=
d\left(
dH-\sum_{i=1}^q dH(R_i)\lambda_i
\right).
\]

Therefore,
\[
L_{X_H}\Omega
=
-\sum_{i=1}^q d(R_i(H))\wedge\lambda_i
-
\sum_{i=1}^qR_i(H)d\lambda_i.
\]
\end{remark}

\begin{example}\label{E33}
Consider the non-uniform $2$-contact manifold $(M=\mathbb R^4,\vec\lambda=(\lambda_1,\lambda_2))$ with coordinates $(x,y,z_1,z_2)$, where
\[
\lambda_1 = dz_1 + x\,dy, \qquad \lambda_2 = dz_2 + c\,x\,dy, \qquad 0 \neq c \neq\pm 1.
\]
The Reeb distribution is $\mathcal R = \operatorname{span}\{\partial_{z_1}, \partial_{z_2}\}$. Moreover, $\xi = \ker\lambda_1 \cap \ker\lambda_2$ is generated by $\partial_x$ and
\[
Y = \partial_y - x\,\partial_{z_1} - c x\,\partial_{z_2}.
\]

We compute $d\lambda_1 = dx \wedge dy$ and $d\lambda_2 = c\,dx \wedge dy$. Hence
\[
\Omega = \frac12(d\lambda_1 + d\lambda_2) = \frac{1+c}{2}\,dx \wedge dy.
\]
Therefore, $\Omega|_\xi$ is non-degenerate whenever $c \neq -1$. Consequently, the structure is admissible but non-uniform.


Let $H = (x^2+y^2+z_1^2+z_2^2)/2$. Then $dH = x\,dx + y\,dy + z_1\,dz_1 + z_2\,dz_2$. The Reeb vector fields are $R_1 = \partial_{z_1}$ and $R_2 = \partial_{z_2}$, hence $dH(R_1)=z_1$ and $dH(R_2)=z_2$. Therefore,
\[
\beta = dH - z_1\lambda_1 - z_2\lambda_2.
\]
Expanding,
\begin{align*}
\beta
&= x\,dx + y\,dy + z_1\,dz_1 + z_2\,dz_2 \\
&\quad - z_1(dz_1 + x\,dy) - z_2(dz_2 + c\,x\,dy) \\
&= x\,dx + \bigl(y - z_1 x - c z_2 x\bigr)\,dy.
\end{align*}

We now determine the horizontal vector field $\widetilde X_H \in \Gamma(\xi)$ satisfying $i_{\widetilde X_H}\Omega = \beta$. Write $\widetilde X_H = A\,\partial_x + B\,Y$. Using $i_{\partial_x}(dx \wedge dy) = dy$ and $i_Y(dx \wedge dy) = -dx$, we obtain
\[
i_{\widetilde X_H}\Omega = \frac{1+c}{2}\bigl(A\,dy - B\,dx\bigr).
\]
Comparing with $\beta$, we find
\[
\frac{1+c}{2}A = y - z_1 x - c z_2 x, \qquad -\frac{1+c}{2}B = x.
\]
Thus
\[
A = \frac{2(y - z_1 x - c z_2 x)}{1+c}, \qquad B = -\frac{2x}{1+c}.
\]
Therefore,
\[
\widetilde X_H = \frac{2(y - z_1 x - c z_2 x)}{1+c}\,\partial_x - \frac{2x}{1+c}\,Y.
\]

The effective Hamiltonian vector field is
\[
X_H = \widetilde X_H - H\bigl(\partial_{z_1} + \partial_{z_2}\bigr).
\]
Substituting the expression of $Y$, we obtain
\begin{align*}
X_H
= \frac{2(y - z_1 x - c z_2 x)}{1+c}\,\partial_x 
- \frac{2x}{1+c}\bigl(\partial_y - x\,\partial_{z_1} - c x\,\partial_{z_2}\bigr) 
 - H\,\partial_{z_1} - H\,\partial_{z_2}.
\end{align*}

Hence the equations of motion are
\begin{align*}
\dot x &= \frac{2(y - z_1 x - c z_2 x)}{1+c}, \\
\dot y &= -\frac{2x}{1+c}, \\
\dot z_1 &= \frac{2x^2}{1+c} - H, \\
\dot z_2 &= \frac{2c x^2}{1+c} - H.
\end{align*}
\end{example}

\subsection{Dynamics on admissible non-uniform $q$-contact structures}

In the uniform theory, the $q$-contact bracket may be defined using any of the contact forms $\lambda_i$, since all associated symplectic structures coincide. In the non-uniform setting, however, the quantities
\[
\lambda_i([X_f^i,X_g^i])
\]
depend on the geometric channel, and there is no canonical way to combine them into a single scalar quantity.

For admissible non-uniform structures, the effective Hamiltonian vector field $X_f$ provides a natural intrinsic dynamics. This motivates the following definition.

\begin{definition}
Let $(M,\vec\lambda)$ be an admissible non-uniform $q$-contact manifold. The \emph{non-uniform $q$-contact bracket} of two smooth functions $f,g\in C^\infty(M)$ is defined by
\begin{equation}\label{QJA}
\{f,g\}
:=
-X_f(g)-g\sum_{i=1}^qR_i(f).
\end{equation}
\end{definition}

\begin{proposition}
The non-uniform $q$-contact bracket is skew-symmetric.
\end{proposition}

\begin{proof}
Since $\Omega = \frac{1}{q}\sum_{i=1}^q d\lambda_i$, the definition of the effective Hamiltonian vector field yields
\[
i_{X_f}\Omega = df - \sum_{i=1}^q R_i(f)\,\lambda_i.
\]
Evaluating this identity on $X_g$, we obtain
\[
\Omega(X_f,X_g) = X_g(f) - \sum_{i=1}^q R_i(f)\,\lambda_i(X_g).
\]
Since $\lambda_i(X_g) = -g$, it follows that
\[
\Omega(X_f,X_g) = X_g(f) + g\sum_{i=1}^q R_i(f).
\]
Similarly,
\[
\Omega(X_g,X_f) = X_f(g) + f\sum_{i=1}^q R_i(g).
\]
Because $\Omega$ is skew-symmetric, $\Omega(X_f,X_g) + \Omega(X_g,X_f) = 0$. Therefore,
\[
X_f(g) + X_g(f) + f\sum_{i=1}^q R_i(g) + g\sum_{i=1}^q R_i(f) = 0.
\]
Using Definition \eqref{QJA}, we conclude that $\{f,g\} + \{g,f\} = 0$.
\end{proof}

\begin{remark}
Unlike the Poisson bracket of symplectic geometry, the bracket \eqref{QJA} is not expected to satisfy either the Jacobi identity or the Leibniz rule in general. Consequently, admissible non-uniform $q$-contact manifolds do not naturally carry a Poisson structure.

It would be interesting to investigate whether this bracket gives rise to a more general algebraic structure, such as a Jacobi or almost Jacobi structure, and whether it induces geometric structures analogous to those arising in Poisson geometry. These questions lie beyond the scope of the present work and will be addressed elsewhere.
\end{remark}

\begin{remark}
In the uniform case, $d\lambda_i = d\lambda_1$ for $i=1,\dots,q$, the non-uniform $q$-contact bracket reduces to the $q$-contact bracket introduced in Section~2.
\end{remark}
Next, we present two properties of the non-uniform $q$-contact bracket.
\begin{proposition}
For every $f,g,h \in C^\infty(M)$, the non-uniform $q$-contact bracket satisfies:
\begin{enumerate}
\item $\{f,gh\} = g\{f,h\} + h\{f,g\} + gh\sum_{i=1}^q R_i(f)$,
\item $\{f,1\} = -\sum_{i=1}^q R_i(f)$.
\end{enumerate}
\end{proposition}

\begin{proof}
Using the definition of the bracket together with the Leibniz rule,
\begin{align*}
\{f,gh\}
&= -X_f(gh) - gh\sum_{i=1}^q R_i(f) \\
&= -gX_f(h) - hX_f(g) - gh\sum_{i=1}^q R_i(f).
\end{align*}
On the other hand, $g\{f,h\} = -gX_f(h) - gh\sum_{i=1}^q R_i(f)$, and similarly $h\{f,g\} = -hX_f(g) - gh\sum_{i=1}^q R_i(f)$. Adding both expressions gives
\[
g\{f,h\} + h\{f,g\} = -gX_f(h) - hX_f(g) - 2gh\sum_{i=1}^q R_i(f).
\]
Hence,
\[
\{f,gh\} = g\{f,h\} + h\{f,g\} + gh\sum_{i=1}^q R_i(f).
\]
Finally, $\{f,1\} = -X_f(1) - \sum_{i=1}^q R_i(f) = -\sum_{i=1}^q R_i(f)$.
\end{proof}
\begin{remark}
The first identity shows that the non-uniform $q$-contact bracket satisfies a modified Leibniz rule. The additional term $gh\sum_{i=1}^q R_i(f)$ measures the deviation from the derivation property characteristic of Poisson brackets. However, we observe that if we restrict to the subspace of functions
\[
C_{\vec\lambda}^\infty(M) = \left\{ f \in C^\infty(M) \;\middle|\; \sum_{i=1}^q R_i(f) = 0 \right\},
\]
then the non-uniform $q$-contact bracket satisfies the derivation property. Consequently, this function subspace, equipped with the bracket, becomes a Leibniz algebra.
\end{remark}

Moreover, by equation \eqref{QJA}, we know that along the Hamiltonian flow generated by $H$,
\[
\frac{dF}{dt}=X_H(F) = -\{H,F\} - F\sum_{i=1}^q R_i(H).
\]
\begin{corollary}\label{Coro}
The evolution of the Hamiltonian function is given by $\frac{dH}{dt} = -H\sum_{i=1}^q R_i(H)$.
\end{corollary}

\begin{proof}
Setting $F=H$ in $X_H(F) = -\{H,F\} - F\sum_{i=1}^q R_i(H)$ and using the skew-symmetry of the bracket yields $\{H,H\}=0$. Therefore, $X_H(H) = -H\sum_{i=1}^q R_i(H)$, which proves the result.
\end{proof}

\begin{definition}\label{DissCons}
A function $f \in C^\infty(M)$ is called:
\begin{enumerate}
\item a \emph{dissipated quantity} if $X_H(f) = -f\sum_{i=1}^q R_i(H)$,
\item a \emph{conserved quantity} if $X_H(f) = 0$.
\end{enumerate}
\end{definition}
\begin{remark}
In general, by Corollary \ref{Coro}, the Hamiltonian function is not conserved along the flow.  However, whenever $\sum_{i=1}^q R_i(H) = 0$, the Hamiltonian function becomes a conserved quantity.
\end{remark}
\begin{remark}
Unlike the uniform case, there is generally no canonical Darboux coordinate system simultaneously normalizing all the forms $\lambda_1,\dots,\lambda_q$. This reflects the presence of several distinct symplectic structures on the horizontal distribution $\xi$.

Nevertheless, one may introduce local coordinates adapted to the splitting $TM = \mathcal R \oplus \xi$. The geometry is then encoded by the family of horizontal symplectic forms $d\lambda_i|_\xi$, or equivalently by the associated structure endomorphisms $B_i: \xi \to \xi$.


In particular, when the structure becomes uniform, i.e. $d\lambda_1 = \cdots = d\lambda_q$, one recovers the Darboux coordinates of uniform $q$-contact geometry and the effective Hamiltonian dynamics reduces to the uniform theory.
\end{remark}
\subsection{Volume evolution and dynamical symmetries}

We now investigate the behavior of the effective Hamiltonian flow with respect to natural volume forms on an admissible non-uniform $q$-contact manifold.

Let $(M,\vec\lambda)$ be an admissible non-uniform $q$-contact manifold, and define $\Omega := \frac{1}{q}\sum_{i=1}^q d\lambda_i$. Since $\Omega|_\xi$ is non-degenerate, the form
\[
\mu := \lambda_1 \wedge \cdots \wedge \lambda_q \wedge \Omega^n
\]
defines a volume form on $M$.

\begin{proposition}\label{PropLieMu}
The Lie derivative of the volume form $\mu$ along the effective Hamiltonian vector field $X_H$ is given by
\begin{align*}
L_{X_H}\mu
&= \sum_{i=1}^q
\lambda_1 \wedge \cdots \wedge (L_{X_H}\lambda_i) \wedge \cdots \wedge \lambda_q \wedge \Omega^n \\
&\quad - n\,\lambda_1 \wedge \cdots \wedge \lambda_q \wedge \left( \sum_{i=1}^q R_i(H)\,d\lambda_i \right) \wedge \Omega^{n-1}.
\end{align*}
\end{proposition}

\begin{proof}
Using the Leibniz rule,
\begin{align*}
L_{X_H}\mu
&= \sum_{i=1}^q
\lambda_1 \wedge \cdots \wedge (L_{X_H}\lambda_i) \wedge \cdots \wedge \lambda_q \wedge \Omega^n \\
&\quad + \lambda_1 \wedge \cdots \wedge \lambda_q \wedge L_{X_H}(\Omega^n).
\end{align*}
Since $L_{X_H}(\Omega^n) = n(L_{X_H}\Omega) \wedge \Omega^{n-1}$, it remains to compute $L_{X_H}\Omega$.

Using Cartan's identity together with $i_{X_H}\Omega = dH - \sum_{i=1}^q R_i(H)\,\lambda_i$, we obtain
\[
L_{X_H}\Omega = d(i_{X_H}\Omega) = -\sum_{i=1}^q d(R_i(H)) \wedge \lambda_i - \sum_{i=1}^q R_i(H)\,d\lambda_i.
\]
After wedging with $\lambda_1 \wedge \cdots \wedge \lambda_q$, all terms containing $d(R_i(H)) \wedge \lambda_i$ vanish identically. Since each term already contains the factor $\lambda_i$ inside $\lambda_1 \wedge \cdots \wedge \lambda_q$, the wedge product with $d(R_i(H)) \wedge \lambda_i$ vanishes identically. Therefore,
\[
\lambda_1 \wedge \cdots \wedge \lambda_q \wedge L_{X_H}(\Omega^n)
=
- n\,\lambda_1 \wedge \cdots \wedge \lambda_q \wedge \left( \sum_{i=1}^q R_i(H)\,d\lambda_i \right) \wedge \Omega^{n-1}.
\]
Substituting this expression into the Leibniz formula yields the result.
\end{proof}
\begin{remark}
Proposition~\ref{PropLieMu} provides the evolution of the natural volume form along the effective Hamiltonian flow. Such evolution formulas are closely related to the problem of the existence of invariant measures for dissipative Hamiltonian systems. In the contact setting, invariant measures have been extensively studied (see, for instance, \cite{BravettiDeLeonMarreroPadron}), where suitable rescalings of the canonical volume form lead to invariant densities. It would be interesting to investigate whether analogous constructions exist for admissible non-uniform $q$-contact Hamiltonian systems. We leave this question for future work.
\end{remark}
\begin{remark}
Unlike the uniform case, the evolution of the volume form receives independent contributions from the different geometric channels $\lambda_1,\dots,\lambda_q$ and $d\lambda_1,\dots,d\lambda_q$. Consequently, the effective Hamiltonian flow does not preserve the volume form $\mu$ in general. Equivalently, the divergence of the effective Hamiltonian vector field with respect to $\mu$ is not expected to vanish.
\end{remark}

\begin{remark}
In the uniform case, $d\lambda_1 = \cdots = d\lambda_q$, all geometric channels coincide and $\Omega = d\lambda_1$, so $\Omega^n = (d\lambda_1)^n$. Therefore, the volume form $\mu = \lambda_1 \wedge \cdots \wedge \lambda_q \wedge \Omega^n$ coincides with the natural volume form of the uniform theory.

By \eqref{LX} and \eqref{LXX}, the evolution of the volume form becomes proportional to the volume form itself. In particular, if $\sum_{i=1}^q R_i(H) = 0$, then the Hamiltonian vector field is divergence-free with respect to $\mu$.
\end{remark}

We now establish the analogue of Noether's theorem for admissible non-uniform $q$-contact structures.

\begin{theorem}[Non-uniform $q$-contact Noether theorem]
For a smooth function $F \in C^\infty(M)$, the following statements are equivalent:
\begin{enumerate}
\item $\{H,F\} = 0$;
\item $F$ is a dissipated quantity, namely $X_H(F) = -F\sum_{i=1}^q R_i(H)$.
\end{enumerate}
Moreover, $F$ is a conserved quantity if and only if $\{H,F\} = -F\sum_{i=1}^q R_i(H)$.
\end{theorem}

\begin{proof}
By definition of the non-uniform $q$-contact bracket,
\[
\{H,F\} = -X_H(F) - F\sum_{i=1}^q R_i(H).
\]
Therefore, $\{H,F\} = 0$ if and only if $X_H(F) = -F\sum_{i=1}^q R_i(H)$. Likewise, $X_H(F) = 0$ if and only if $\{H,F\} = -F\sum_{i=1}^q R_i(H)$.
\end{proof}
\subsubsection{Effective structure endomorphism}

Recall that the horizontal distribution $\xi$ is endowed with the reference symplectic form $\omega = d\lambda_1|_\xi$. For each $i=1,\dots,q$, the corresponding symplectic structure $d\lambda_i|_\xi$ is related to $\omega$ through the structure endomorphism $B_i: \xi \to \xi$, defined by
\[
d\lambda_i(U,V) = \omega(B_i U, V), \qquad U,V \in \Gamma(\xi).
\]
The effective symplectic form $\Omega = \frac{1}{q}\sum_{i=1}^q d\lambda_i$ therefore satisfies
\[
\Omega(U,V) = \omega\!\left( \frac1q \sum_{i=1}^q B_i\,U, V \right), \qquad U,V \in \Gamma(\xi).
\]
Consequently, the admissible non-uniform geometry is governed by the effective structure endomorphism $B = \frac1q \sum_{i=1}^q B_i$.

\begin{remark}
The admissibility condition is equivalent to the invertibility of the effective structure endomorphism $B$. Indeed, $\Omega(U,V) = \omega(BU,V)$, and since $\omega$ is non-degenerate on $\xi$, the restriction $\Omega|_\xi$ is non-degenerate if and only if $B$ is invertible.
\end{remark}

\begin{remark}
In the uniform case, $B_i = \operatorname{Id}_\xi$ for $i=1,\dots,q$, and therefore $B = \frac1q \sum_{i=1}^q B_i = \operatorname{Id}_\xi$. Consequently, $\Omega|_\xi = \omega$, recovering the effective symplectic structure of the uniform theory.
\end{remark}

\begin{remark}
The effective Hamiltonian vector field $X_H$ is determined by
\[
i_{X_H}\Omega = dH - \sum_{i=1}^q dH(R_i)\,\lambda_i,
\]
together with $\lambda_i(X_H) = -H$ for $i=1,\dots,q$. Restricted to the horizontal distribution $\xi$, the dynamics is controlled by the effective operator $B = \frac1q \sum_{i=1}^q B_i$. Deviations of $B$ from the identity encode the anisotropy of the non-uniform $q$-contact structure and determine how the different geometric channels deform the Hamiltonian evolution.

In particular, if $\frac1q \sum_{i=1}^q B_i = c\,\operatorname{Id}_\xi$ with $c \neq 0$, then $\Omega|_\xi = c\,\omega$. In this case, the horizontal part of the effective Hamiltonian vector field is proportional to the horizontal component associated with the first geometric channel. More precisely, $Y_H = \frac1c\,Y_H^1$, where $Y_H, Y_H^1 \in \Gamma(\xi)$ denote the corresponding horizontal components.
\end{remark}
\subsection{Towards a More General Structure}
In fact, the $q$-contact structure can be subsumed into a more general framework. In the following, we present some discussions on this general structure.
	\begin{definition}(q-general manifolds)
		Let $n,q$ be positive integers and consider a $2n+q$-dimensional differentiable manifold $M$. A $q$-general structure on $M$ consists of a collection $\vec{\lambda}=(\lambda_1,\dots,\lambda_q)$ of $q$ (pointwise) linearly independent non-vanishing $1$-forms $\lambda_i$, $i=1,\dots,q$, and a collection $\vec{\omega}=(\Omega_1,\dots,\Omega_q)$ of $n$ non-vanishing $2$-forms $\Omega_j$, $j=1,\dots,n$, together with a splitting 
		$$TM = \mathcal{R} \oplus \xi$$ 
		of the tangent bundle, satisfying the following conditions:
		
		(i) $\xi := \bigcap_{i=1}^q \ker \lambda_i$;
		
		(ii) for every $j$, the restriction $\Omega_j|_\xi$ is non-degenerate;
		
		(iii) for every $j$, one has $\ker \Omega_j = \mathcal{R}$.
	\end{definition}
	
	\begin{remark}
		The linear independence of the $\lambda_i$ implies that $\xi$ has constant rank $2n$, and condition (ii) is equivalent to $\Omega_1 \wedge \cdots \wedge \Omega_n|_\xi \neq 0$; in other words, $(\xi, \Omega_j)$, $j=1,\dots,n$, are symplectic bundles over $M$.
	\end{remark}
	\begin{remark}
		In fact, we can see that $\Omega_{i_1}\wedge\cdots\wedge\Omega_{i_n}|_\xi\neq 0,i_j\in\{1,...,n\},j=1,...,n.$
	\end{remark}
	A manifold endowed with such structure is called a $q$-$general\;manifold$ and denoted by $(M,\vec{\lambda}\oplus \vec{\omega},\mathcal R\oplus\xi)$. We call the collections $\{\lambda_i \}$ and $\{\Omega_j\}$  $adapted\;coframes$ for the $q$-general structure. The $q$-form
$$\lambda:=\lambda_1\wedge\cdots\wedge\lambda_q\neq0$$
is called the $characteristic\;form.$ The $2n$-form
$$\Omega_1\wedge\cdots\wedge\Omega_n$$
is called the $cocharacteristic\;form.$
The bundles $\mathcal R$ and $\xi$ are called the $Reeb\;distribution$ and $q$-$general\;distribution,$ respectively. The elements of $\xi$ will be called $horizontal\;vector\;field.$
\begin{remark}
	When $\Omega := \Omega_1 = \cdots = \Omega_n$, in \cite{de1} the authors called the $(q+1)$-tuple $(\Omega, \lambda_1, \dots, \lambda_q)$ an almost cosymplectic structure of order $q$. Moreover, if $d\Omega_1 = 0$, the authors called the $(q+1)$-tuple $(\Omega, \lambda_1, \dots, \lambda_q)$ a partially cosymplectic structure of order $q$. Furthermore, if $\lambda_1, \dots, \lambda_q$ are also closed, then $(\Omega, \lambda_1, \dots, \lambda_q)$ is called cosymplectic of order $q$, while in \cite{Leok2} such a structure is called a $q$-cosymplectic structure (which is the terminology adopted in this paper). When $\Omega_1 = \cdots = \Omega_n = d\lambda_1 = \cdots = d\lambda_q$, we return to the case of uniform $q$-contact structure. When $\Omega_i=d\lambda_j,i\in\{1,...,n\},j\in\{1,...,q\}$, we return to the case of general $q$-contact manifold which we will discuss later.
\end{remark}

\begin{proposition}
	Given a $q$-general manifold $(M,\vec{\lambda}\oplus \vec{\omega},\mathcal R\oplus\xi)$, for any $1\leq j\leq n,$ the vector bundle morphisms given by
	\[
	\flat^j : TM \longrightarrow T^*M, \quad j=1,\dots,n
	\]
	\[
	X \longmapsto i_X \Omega_j + \sum_{k=1}^q \lambda_k(X) \lambda_k,
	\]
	are isomorphisms. Consequently, they induce isomorphisms of \( C^\infty(M) \)-modules between the vector fields on \( M \) and the $1$-forms on \( M \).
\end{proposition}

\begin{proof}
	It suffices to prove that for any \(1 \leq j \leq n\), the map \( \flat^j_x : T_xM \longrightarrow T_x^*M \) is an isomorphism for each \( x \in M \).
	
	Since \(T_x^* M\) and \(T_x M\) have the same dimension, it suffices to show that \(\flat^j\) is injective.  
	Suppose there exists a non-zero tangent vector \(X \in T_x M\) such that \(\flat^j(X) = 0\).  
	Then \(\flat^j(X)(X) = \sum_{k=1}^q (\lambda_k(X))^2 = 0\), so we conclude that \(\lambda_k(X) = 0\) for all \(k\), which means \(X \in \xi\). Moreover, we also have \(i_X \Omega_j = 0\), thus \(X = 0\).  
	
	Conversely, by definition, we can decompose \(T^*M = \tilde\xi \oplus \tilde{\mathcal R}\) such that any element of \(\tilde\xi\) pairs trivially with elements of \(\mathcal R\), and any element of \(\tilde{\mathcal R}\) pairs trivially with elements of \(\xi\). Here, \(\Omega_j\) provides an isomorphism between \(\xi\) and \(\tilde\xi\).  
	Thus, for any \(\eta \in \Lambda^1(TM)\), we have \(\eta - \sum_{k=1}^q R_k(\eta) \lambda_k \in \tilde\xi\). Using \(\Omega_j\), we obtain a vector field \(\tilde X \in \xi\). Letting \(X = \tilde X + \sum_{k=1}^q R_k(\eta) R_k\), we can verify that \(\flat^j(X) = \eta\).  
	Hence \(\flat^j\) is surjective and therefore an isomorphism.
\end{proof}
\begin{corollary}
	Given a $q$-general manifold $(M,\vec{\lambda}\oplus \vec{\omega},\mathcal R\oplus\xi)$, for any $1\leq j\leq n$, there exist a unique vector fields \(R^j_k\), \(k = 1, \ldots, p\), such that:
	\[
	i_{R_k^j} \Omega_j = 0, \qquad i_{R_k^j} \eta_l = \delta_{kl}.
	\]
	The family \(\{R_k^j \mid k = 1, \ldots, p\}\) will be called the family of Reeb vector fields about $\Omega_j$ of the structure.
\end{corollary}

\begin{proof}
	It suffices to set \(R_k^j = (\flat^j)^{-1}(\eta_k)\).
\end{proof}
\begin{remark}
	In \cite{de1}, the authors proved that, given an almost cosymplectic structure of order \( q \) \( (\Omega, \lambda_1, \cdots, \lambda_q) \) on \( M \), there exist unique vector fields \( {R}_k, k = 1, \ldots, q \) such that:
	\[
	i_{{R}_k} \omega = 0, \quad i_{{R}_k} \eta_j = \delta_{kj}.
	\]
	The family \( \{ {R}_k, k = 1, \ldots, q \} \) is called the family of Reeb vector fields of the structure. In \cite{Leok2}, the authors proved that for a \(q\)-cosymplectic manifold, the family \( \{ {R}_k, k = 1, \ldots, q \} \) further satisfies \( [R_i,R_j]=0 \). The same conclusion can also be obtained for uniform \(q\)-contact manifolds (see \cite{Finamore}).
\end{remark}
\section{Application: An admissible non-uniform dissipative system with two internal channels}

Many dissipative materials exhibit internal mechanisms that influence the macroscopic response of the system. Typical examples arise in viscoelasticity, rheology, and continuum thermodynamics with internal variables, where part of the stored energy is associated with microscopic processes that are not directly observable at the macroscopic level.

A common modeling strategy consists of introducing internal variables describing different dissipation channels. These variables account for the presence of internal microstructural mechanisms such as relaxation processes, rearrangements of the material structure, or energy exchange between different scales.

We consider a mechanical system with two internal channels. Let
$\varepsilon$ denote the macroscopic strain and let $p$ be its conjugate momentum. In addition, we introduce two internal variables $\alpha_1,\alpha_2,$ representing two distinct internal modes. Finally, we introduce two dissipative variables
$s_1, s_2,$ associated with the irreversible contribution of each channel.

The phase space is therefore

\[
M=\mathbb R^6
(\varepsilon,p,\alpha_1,\alpha_2,s_1,s_2).
\]

The energy of the system is described by the following Hamiltonian:

\[
H
=
\frac{p^2}{2m}
+
\frac{k}{2}
(\varepsilon-\alpha_1-\alpha_2)^2
+
\frac{k_1}{2}\alpha_1^2
+
\frac{k_2}{2}\alpha_2^2
+
\kappa\alpha_1\alpha_2
+
\mu_1 s_1
+
\mu_2 s_2.
\]
\noindent
where \(m>0\) is an effective mass parameter, \(k>0\) is the elastic stiffness associated with the observable strain, \(k_1>0\) and \(k_2>0\) measure the stiffness of the two internal modes, \(\kappa\in\mathbb R\) describes the coupling between the internal channels, \(\mu_1>0\) and \(\mu_2>0\) are dissipation coefficients associated with the two dissipative variables \(s_1\) and \(s_2\).

The quantity $\varepsilon-\alpha_1-\alpha_2$ represents the effective elastic deformation stored by the system. Consequently, the total strain is decomposed into an observable contribution and two internal components. The variables \(\alpha_1\) and \(\alpha_2\) therefore act as hidden degrees of freedom capable of storing and exchanging energy with the observable mechanical mode. The terms $\mu_1s_1,
\mu_2s_2,$ introduce dissipation through the two geometric channels of the non-uniform \(2\)-contact structure.

We now show that this system naturally determines an admissible non-uniform $2$-contact Hamiltonian system.

Define the one-forms:
\begin{equation*}
\lambda_1 = ds_1 + p\,d\varepsilon + \alpha_1\,d\alpha_2,\qquad
\lambda_2 = ds_2 + c_1 p\,d\varepsilon + c_2 \alpha_1\,d\alpha_2,
\end{equation*}
where $c_1,c_2\in\mathbb R\setminus\{0,-1\}$, $(c_1,c_2)\neq(1,1)$. The coefficients $c_1,c_2$ quantify the geometric asymmetry between the two channels; $(c_1,c_2)=(1,1)$ corresponds to the uniform case, otherwise non-uniform.

The exterior derivatives are
\begin{equation*}
d\lambda_1 = dp\wedge d\varepsilon + d\alpha_1\wedge d\alpha_2,\qquad
d\lambda_2 = c_1\,dp\wedge d\varepsilon + c_2\,d\alpha_1\wedge d\alpha_2.
\end{equation*}
Since $(c_1,c_2)\neq(1,1)$, $d\lambda_1\neq d\lambda_2$, hence the structure is genuinely non-uniform.

Consider a vector field
\begin{equation*}
X = a\,\partial_\varepsilon + b\,\partial_p + u\,\partial_{\alpha_1} + v\,\partial_{\alpha_2} + r_1\,\partial_{s_1} + r_2\,\partial_{s_2}.
\end{equation*}
Then
\begin{equation*}
i_X d\lambda_1 = b\,d\varepsilon - a\,dp + u\,d\alpha_2 - v\,d\alpha_1.
\end{equation*}
Thus $i_X d\lambda_1=0$ iff $a=b=u=v=0$. Therefore
\begin{equation*}
\ker d\lambda_1 = \operatorname{span}\{\partial_{s_1}, \partial_{s_2}\}.
\end{equation*}
The same computation gives $\ker d\lambda_2 = \operatorname{span}\{\partial_{s_1}, \partial_{s_2}\}$. Consequently,
\begin{equation*}
\mathcal R = \operatorname{span}\{R_1,R_2\},\qquad R_1=\partial_{s_1},\quad R_2=\partial_{s_2}.
\end{equation*}

The horizontal distribution is $\xi = \ker\lambda_1\cap\ker\lambda_2$. The conditions $\lambda_1(X)=0,\ \lambda_2(X)=0$ are equivalent to
\begin{equation*}
r_1 = -p a - \alpha_1 v,\qquad r_2 = -c_1 p a - c_2 \alpha_1 v.
\end{equation*}
Thus every vector field in $\xi$ is uniquely determined by $(a,b,u,v)$, hence $\dim\xi=4$.

A convenient basis of $\xi$ is
\begin{align*}
E_1 &= \partial_\varepsilon - p\,\partial_{s_1} - c_1 p\,\partial_{s_2}, &
E_2 &= \partial_p, \\
E_3 &= \partial_{\alpha_1}, &
E_4 &= \partial_{\alpha_2} - \alpha_1\partial_{s_1} - c_2\alpha_1\partial_{s_2}.
\end{align*}

 With respect to this basis, the restriction \(d\lambda_1|_\xi\) is represented by the matrix
\[
\begin{pmatrix}
0&-1&0&0\\
1&0&0&0\\
0&0&0&1\\
0&0&-1&0
\end{pmatrix},
\]
whose determinant is equal to one. Hence $d\lambda_1|_\xi$ is non-degenerate. Similarly, $d\lambda_2|_\xi$ is represented by

\[
\begin{pmatrix}
0&-c_1&0&0\\
c_1&0&0&0\\
0&0&0&c_2\\
0&0&-c_2&0
\end{pmatrix},
\]
\noindent
whose determinant is $c_1^2c_2^2.$ Since $c_1\neq0, c_2\neq0,$ the restriction
$d\lambda_2|_\xi$ is also non-degenerate. Therefore both forms induce symplectic structures on the horizontal distribution. Taking $\omega=d\lambda_1|_\xi$
as the reference symplectic form, the structure endomorphism
$B_2:\xi\to\xi$ is determined by
\[
d\lambda_2(U,V)
=
\omega(B_2U,V),
\qquad
U,V\in\xi.
\]

A direct computation gives
\[
B_2
=
\operatorname{diag}
(c_1,c_1,c_2,c_2).
\]
The effective structure endomorphism is therefore
\[
B
=
\frac12(I+B_2),
\]
namely
\[
B
=
\operatorname{diag}
\left(
\frac{1+c_1}{2},
\frac{1+c_1}{2},
\frac{1+c_2}{2},
\frac{1+c_2}{2}
\right).
\]
Consequently,
\[
\det(B)
=
\left(\frac{1+c_1}{2}\right)^2
\left(\frac{1+c_2}{2}\right)^2.
\]
Thus \(B\) is invertible if and only if
\[
c_1\neq-1,
\qquad
c_2\neq-1.
\]
The effective two-form is
\[
\Omega
=
\frac12(d\lambda_1+d\lambda_2)
=
\frac12
\Big(
(1+c_1)\,dp\wedge d\varepsilon
+
(1+c_2)\,d\alpha_1\wedge d\alpha_2
\Big).
\]
Its restriction to \(\xi\) is represented by
\[
\begin{pmatrix}
0&-\frac{1+c_1}{2}&0&0\\
\frac{1+c_1}{2}&0&0&0\\
0&0&0&\frac{1+c_2}{2}\\
0&0&-\frac{1+c_2}{2}&0
\end{pmatrix},
\]
whose determinant equals
\[
\frac{(1+c_1)^2(1+c_2)^2}{16}.
\]
Therefore $\Omega|_\xi$ is non-degenerate precisely  when

\[
c_1\neq-1,
\qquad
c_2\neq-1.
\]

It follows that the non-uniform $2$-contact structure
\[
(M,\vec{\lambda},\mathcal R\oplus\xi),\qquad \vec{\lambda}=(\lambda_1,\lambda_2),
\]
with
\[
\mathcal R = \operatorname{span}\left\{\frac{\partial}{\partial s_1},\frac{\partial}{\partial s_2}\right\},
\quad
\xi = \ker\lambda_1 \cap \ker\lambda_2,
\]
is admissible.

Recall that $\mathcal R$ is the Reeb distribution and $TM=\mathcal R\oplus\xi$. Moreover, since $R_1(H)=\mu_1,\ R_2(H)=\mu_2$, Corollary~\ref{Coro} implies that every effective Hamiltonian trajectory satisfies
\[
\frac{dH}{dt} = -(\mu_1+\mu_2)H.
\]
Consequently,
\[
H(t) = H(0)e^{-(\mu_1+\mu_2)t}.
\]
Therefore the total dissipation rate is determined by the sum of the contributions of the two Reeb channels. In particular, the dissipation law is independent of the anisotropy parameters $c_1,c_2$, which only affect the geometry of the admissible structure.
\subsection{Effective Hamiltonian dynamics}

The admissibility of the non-uniform $2$-contact structure guarantees the existence of an effective Hamiltonian vector field associated with $H$. According to Theorem~\ref{TT}, the effective Hamiltonian vector field $X_H$ is uniquely determined by
\[
i_{Y_H}\Omega = dH - \sum_{i=1}^{2}R_i(H)\lambda_i,
\]
where $Y_H\in\xi$ denotes the horizontal component of $X_H$, together with $X_H = Y_H - H(R_1+R_2)$. For this example, $R_1(H)=\mu_1,\ R_2(H)=\mu_2$, hence
\[
i_{Y_H}\Omega = dH - \mu_1\lambda_1 - \mu_2\lambda_2.
\]

Introducing $E=\varepsilon-\alpha_1-\alpha_2$, one obtains
\begin{align*}
dH &=
kE\,d\varepsilon + \frac{p}{m}\,dp
+ (-kE+k_1\alpha_1+\kappa\alpha_2)\,d\alpha_1 \\
&\quad + (-kE+k_2\alpha_2+\kappa\alpha_1)\,d\alpha_2
+ \mu_1\,ds_1 + \mu_2\,ds_2.
\end{align*}
Hence
\begin{align*}
\beta:=dH-\mu_1\lambda_1-\mu_2\lambda_2
&=
\bigl(kE-(\mu_1+c_1\mu_2)p\bigr)\,d\varepsilon
+ \frac{p}{m}\,dp \\
&\quad + (-kE+k_1\alpha_1+\kappa\alpha_2)\,d\alpha_1
+ \bigl(-kE+k_2\alpha_2+\kappa\alpha_1-(\mu_1+c_2\mu_2)\alpha_1\bigr)\,d\alpha_2.
\end{align*}

Let $Y_H = A E_1 + B E_2 + U E_3 + V E_4$. Evaluating $\beta$ on $(E_1,E_2,E_3,E_4)$ and solving $i_{Y_H}\Omega=\beta$ yields
\begin{align*}
A &= -\frac{2}{1+c_1}\frac{p}{m}, &
B &= \frac{2}{1+c_1}\bigl(kE-(\mu_1+c_1\mu_2)p\bigr), \\
U &= \frac{2}{1+c_2}\bigl(-kE+k_2\alpha_2+\kappa\alpha_1-(\mu_1+c_2\mu_2)\alpha_1\bigr), &
V &= \frac{2}{1+c_2}\bigl(kE-k_1\alpha_1-\kappa\alpha_2\bigr).
\end{align*}
Therefore
\[
X_H = A E_1 + B E_2 + U E_3 + V E_4 - H\,\partial_{s_1} - H\,\partial_{s_2}.
\]

Using the explicit basis vectors
\[
E_1=\partial_\varepsilon - p\,\partial_{s_1} - c_1p\,\partial_{s_2},\quad
E_2=\partial_p,\quad
E_3=\partial_{\alpha_1},\quad
E_4=\partial_{\alpha_2} - \alpha_1\partial_{s_1} - c_2\alpha_1\partial_{s_2},
\]
the effective Hamiltonian vector field can be written as
\[
X_H = \dot\varepsilon\,\partial_\varepsilon + \dot p\,\partial_p + \dot\alpha_1\,\partial_{\alpha_1} + \dot\alpha_2\,\partial_{\alpha_2} + \dot s_1\,\partial_{s_1} + \dot s_2\,\partial_{s_2}.
\]
Consequently, the effective dynamics is governed by
\begin{align*}
\dot\varepsilon &= -\frac{2}{1+c_1}\frac{p}{m}, \\
\dot p &= \frac{2}{1+c_1}\bigl(k(\varepsilon-\alpha_1-\alpha_2)-(\mu_1+c_1\mu_2)p\bigr), \\
\dot\alpha_1 &= \frac{2}{1+c_2}\bigl(-k(\varepsilon-\alpha_1-\alpha_2)+k_2\alpha_2+\kappa\alpha_1-(\mu_1+c_2\mu_2)\alpha_1\bigr), \\
\dot\alpha_2 &= \frac{2}{1+c_2}\bigl(k(\varepsilon-\alpha_1-\alpha_2)-k_1\alpha_1-\kappa\alpha_2\bigr), \\
\dot s_1 &= -pA - \alpha_1 V - H, \\
\dot s_2 &= -c_1 p A - c_2\alpha_1 V - H.
\end{align*}
Substituting $A$ and $V$, we obtain
\begin{align*}
\dot s_1 &=
\frac{2p^2}{m(1+c_1)}
- \frac{2\alpha_1}{1+c_2}\bigl(k(\varepsilon-\alpha_1-\alpha_2)-k_1\alpha_1-\kappa\alpha_2\bigr) - H, \\
\dot s_2 &=
\frac{2c_1p^2}{m(1+c_1)}
- \frac{2c_2\alpha_1}{1+c_2}\bigl(k(\varepsilon-\alpha_1-\alpha_2)-k_1\alpha_1-\kappa\alpha_2\bigr) - H.
\end{align*}

These equations define the effective evolution associated with the admissible non-uniform $2$-contact Hamiltonian structure. In the next section they will be integrated numerically to investigate the dissipative behaviour predicted by the theory and the interaction between the mechanical and dissipative sectors.

\subsection{Numerical experiments}
The effective Hamiltonian dynamics derived in the previous subsection predicts a precise dissipation law for the total Hamiltonian,
\[
\frac{dH}{dt} = -(\mu_1+\mu_2)H,
\]
hence
\[
H(t) = H(0)e^{-(\mu_1+\mu_2)t}.
\]
Moreover, the rescaled quantity $I(t)=H(t)e^{(\mu_1+\mu_2)t}$ satisfies
\[
\frac{dI}{dt} = e^{(\mu_1+\mu_2)t}\left(\frac{dH}{dt}+(\mu_1+\mu_2)H\right) = 0,
\]
and is therefore conserved along the effective flow.

The purpose of the numerical experiments is twofold. First, we verify the exponential dissipation law satisfied by $H$. Second, we investigate the conservation of $I(t)$, which follows directly from the geometric theory.

The numerical experiments are performed by integrating the effective Hamiltonian system derived above using a fourth-order Runge--Kutta method. The parameters used throughout the simulations are
\[
m=1,\quad k=1,\quad k_1=0.8,\quad k_2=1.2,\quad
\kappa=0.3,\quad \mu_1=0.15,\quad \mu_2=0.25,\quad c_1=2,\quad c_2=0.5.
\]
The initial condition is chosen as
\[
(\varepsilon,p,\alpha_1,\alpha_2,s_1,s_2) = (1,\,0,\,0.2,\,-0.1,\,0,\,0).
\]
\subsubsection{Experiment 1: Verification of the dissipation law}

According to the effective Hamiltonian theory, every trajectory satisfies

\[
\frac{dH}{dt}
=
-(\mu_1+\mu_2)H.
\]

The solution of this equation is

\[
H_{\mathrm{th}}(t)
=
H(0)e^{-(\mu_1+\mu_2)t}.
\]

To verify this prediction, we compare the numerical Hamiltonian obtained from the integration of the effective Hamiltonian system with the theoretical exponential decay.

\begin{figure}[H]
\centering
\includegraphics[width=.75\textwidth]{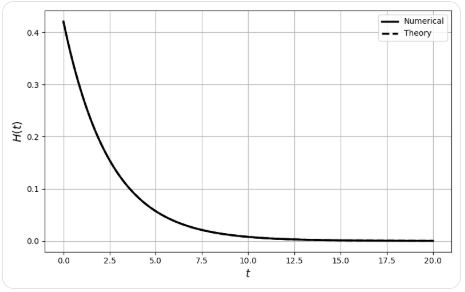}
\vspace{0.4cm}
\includegraphics[width=.75\textwidth]{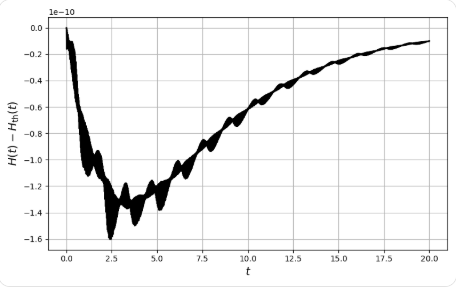}
\caption{Top: comparison between the numerical Hamiltonian \(H(t)\) and the theoretical prediction
\(
H_{\mathrm{th}}(t)=H(0)e^{-(\mu_1+\mu_2)t}.
\)
Bottom: error
\(
E_H(t)=H(t)-H_{\mathrm{th}}(t).
\)
The error remains at machine precision throughout the simulation.}
\label{fig:dissipation}
\end{figure}

Figure~\ref{fig:dissipation} shows the numerical evolution of the Hamiltonian together with the theoretical prediction \(H_{\mathrm{th}}(t)\). The two curves are visually indistinguishable, confirming the dissipation law predicted by Corollary~\ref{Coro}. In order to quantify the agreement, we define the error

\[
E_H(t)
=
H(t)-H_{\mathrm{th}}(t).
\]

The numerical computations show that \(E_H(t)\) remains close to machine precision throughout the simulation interval, providing strong evidence for the validity of the theoretical dissipation law.
The top part of Figure~\ref{fig:dissipation} shows the numerical Hamiltonian together with the theoretical prediction

\[
H_{\mathrm{th}}(t)
=
H(0)e^{-(\mu_1+\mu_2)t}.
\]

The two curves are visually indistinguishable throughout the integration interval. To quantify the agreement, we compute the difference

\[
E_H(t)
=
H(t)-H_{\mathrm{th}}(t).
\]

The bottom part of Figure \ref{fig:dissipation} shows that the error remains of order \(10^{-10}\), which is consistent with the numerical accuracy of the integration scheme. This provides a strong numerical confirmation of the effective dissipation law predicted by the theory.

\subsubsection{Experiment 2: Conservation of the rescaled energy}

The dissipation law established in Corollary~\ref{Coro} implies that the quantity

\[
I(t)
=
H(t)e^{(\mu_1+\mu_2)t}
\]

is conserved along every effective Hamiltonian trajectory. Indeed,

\[
\frac{dI}{dt}
=
e^{(\mu_1+\mu_2)t}
\left(
\frac{dH}{dt}
+
(\mu_1+\mu_2)H
\right),
\]

and therefore

\[
\frac{dI}{dt}=0.
\]

This quantity may be interpreted as the Hamiltonian after removing the universal exponential dissipation generated by the two Reeb channels. To verify this conservation law numerically, we compute \(I(t)\) along the trajectories obtained from the effective Hamiltonian system and compare it with its initial value \(I(0)\).

In addition, we evaluate the conservation error

\[
E_I(t)
=
I(t)-I(0).
\]

If the theoretical prediction is correct, \(I(t)\) should remain constant and \(E_I(t)\) should stay close to machine precision throughout the simulation.

\begin{figure}[H]
\centering
\includegraphics[width=.75\textwidth]{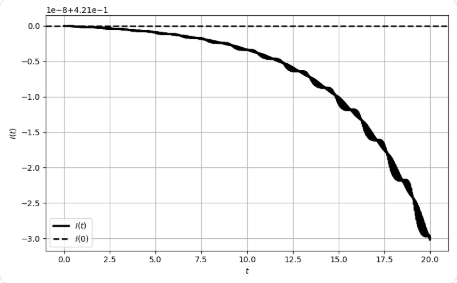}

\vspace{0.4cm}

\includegraphics[width=.75\textwidth]{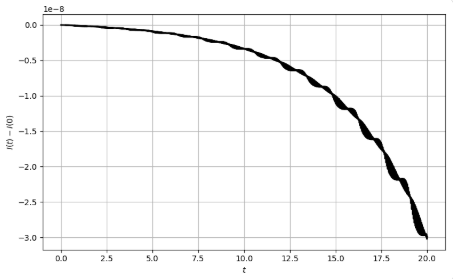}

\caption{Top: evolution of the rescaled quantity
\(
I(t)=H(t)e^{(\mu_1+\mu_2)t}.
\)
The dashed line represents the initial value \(I(0)\).
Bottom: conservation error
\(
E_I(t)=I(t)-I(0).
\)
The variation remains extremely small throughout the simulation interval.}
\label{fig:conserved}
\end{figure}

Figure~\ref{fig:conserved} illustrates the behaviour of the rescaled quantity

\[
I(t)
=
H(t)e^{(\mu_1+\mu_2)t}.
\]

The upper panel shows the numerical evolution of \(I(t)\) together with its initial value \(I(0)\). The two curves are visually indistinguishable over the entire integration interval.

To quantify the conservation error, we compute

\[
E_I(t)
=
I(t)-I(0).
\]

The lower panel shows that the variation of \(I(t)\) remains of order \(10^{-8}\). Since the magnitude of \(I(t)\) is approximately \(4.2\times10^{-1}\), the corresponding relative error is of order \(10^{-7}\).

Therefore, the numerical results confirm that \(I(t)\) remains constant up to numerical accuracy. Together with the exponential decay of the Hamiltonian established in Experiment~1, this illustrates the coexistence of dissipated and conserved quantities predicted by the effective non-uniform \(2\)-contact Hamiltonian dynamics.

\section{Conclusions}

In this work, we introduced the notion of admissible non-uniform $q$-contact structures and developed the corresponding effective Hamiltonian formalism. The construction provides a natural geometric framework for describing systems possessing several dissipation channels with different geometric characteristics.

The key ingredient of the theory is the effective two-form obtained from the family of contact structures. Under the admissibility condition, this form induces a well-defined effective Hamiltonian vector field whose dynamics combines conservative and dissipative effects in a unified geometric setting.

Several fundamental properties of the effective dynamics were established. In particular, we derived explicit evolution equations for the Hamiltonian function and obtained general conditions characterizing dissipated and conserved quantities. These results show that the geometry itself determines how energy is redistributed and dissipated through the different contact channels.

To illustrate the theory, we constructed an admissible non-uniform $2$-contact Hamiltonian system with two internal dissipation mechanisms. The example models a mechanical system with internal variables and provides a concrete realization of the abstract geometric framework. The corresponding effective Hamiltonian vector field was computed explicitly and integrated numerically.

The numerical experiments confirm the theoretical predictions. The Hamiltonian function follows the dissipation law prescribed by the theory, while the rescaled quantity associated with the effective dynamics remains constant up to numerical accuracy. These simulations illustrate the coexistence of dissipated and conserved quantities predicted by the effective non-uniform contact geometry.

The framework developed here opens several directions for future research. Possible developments include the study of non-uniform $q$-contact systems with constraints, symmetry reduction procedures, geometric integrators adapted to the effective dynamics, and applications to continuum mechanics, thermodynamics, and multiphysics systems with several interacting dissipation mechanisms.

Overall, the results show that admissible non-uniform $q$-contact geometry provides a natural extension of contact Hamiltonian mechanics capable of encoding multiple non-equivalent dissipation channels within a single geometric structure.
\section*{Acknowledgments}
We thank the support by NSFC (Grant No.
	12401234). We also acknowledge the financial support of the Ministerio de Ciencia, Innovación y Universidades
(Spain), grant PID2022-125515NB-C21; and the Severo
Ochoa Programme for Centers of Excellence (Grant CEX2023-001347-S).

	$\\$
	
	\noindent$\mathbf{Conflict\;of\;interest\;statement.}$ On behalf of all authors, the corresponding author states that there is no conflict of interest.
	
	$\\$
	\noindent$\mathbf{Data\;availability.}$ Data sharing is not applicable to this article as no new data were created or analyzed in this study.
	
\end{document}